\documentclass[11pt]{article}
\usepackage{graphicx}
\usepackage{amsmath}
\usepackage{amsfonts}
\usepackage{enumerate}
\usepackage{latexsym}
\usepackage{epsfig}

\newtheorem{theorem}{Theorem}[section]
\newtheorem{lemma}[theorem]{Lemma}
\newtheorem{coroll}[theorem]{Corollary}

\def\proofbox{\begin{picture}(6.5,6.5)
\put(0,0){\framebox(6.5,6.5){}}\end{picture}}
\newenvironment{proof}{\noindent{\it Proof.\quad}}{\hfill\proofbox}

\begin{document}
\begin{center} {\bf \Large Complexes of Nonseparating Curves
and}\\

\vspace{0.09in}

{\bf \Large Mapping Class Groups}

\vspace{0.15in}
\large {Elmas Irmak} \\
\vspace{0.09in}
\end{center}



\begin{abstract}Let $R$ be a compact, connected, orientable surface of genus $g$,
$Mod_R^*$ be the extended mapping class group of $R$,
$\mathcal{C}(R)$ be the complex of curves on $R$, and
$\mathcal{N}(R)$ be the complex of nonseparating curves on $R$. We
prove that if $g \geq 2$ and $R$ has at most $g-1$ boundary
components, then a simplicial map $\lambda: \mathcal{N}(R)
\rightarrow \mathcal{N}(R)$ is superinjective if and only if it is
induced by a homeomorphism of $R$. We prove that if $g \geq 2$ and
$R$ is not a closed surface of genus two then
$Aut(\mathcal{N}(R))= Mod_R^*$, and if $R$ is a closed surface of
genus two then $Aut(\mathcal{N}(R))= Mod_R ^* /\mathcal{C}(Mod_R
^*)$. We also prove that if $g=2$ and $R$ has at most one boundary
component, then a simplicial map $\lambda: \mathcal{C}(R)
\rightarrow \mathcal{C}(R)$ is superinjective if and only if it is
induced by a homeomorphism of $R$. As a corollary we prove some
new results about injective homomorphisms from finite index
subgroups to $Mod_R^*$. The last two results complete the author's
previous results to connected orientable surfaces of genus at
least two.
\end{abstract}




\section{Introduction}

Let $R$ be a compact, connected, orientable surface of genus $g$
with $p$ boundary components. The extended mapping class group,
$Mod_R^*$, of $R$ is the group of isotopy classes of all
(including orientation reversing) homeomorphisms of $R$. Let
$\mathcal{A}$ denote the set of isotopy classes of nontrivial
simple closed curves on $R$. The \textit{complex of curves},
$\mathcal{C}(R)$, on $R$ is an abstract simplicial complex,
introduced by Harvey \cite{Har}, with vertex set $\mathcal{A}$
such that a set of $n$ vertices $\{{ \alpha_{1}}, {\alpha_{2}},
..., {\alpha_{n}}\}$ forms an $n-1$ simplex if and only if
${\alpha_{1}}, {\alpha_{2}},..., {\alpha_{n}}$ have pairwise
disjoint representatives. Let $\mathcal{B}$ denote the set of
isotopy classes of nonseparating simple closed curves on $R$. The
\textit{complex of nonseparating curves}, $\mathcal{N}(R)$, is the
subcomplex of $\mathcal{C}(R)$ with the vertex set $\mathcal{B}$
such that a set of $n$ vertices $\{{ \beta_{1}}, {\beta_{2}}, ...,
{\beta_{n}}\}$ forms an $n-1$ simplex if and only if ${\beta_{1}},
{\beta_{2}},..., {\beta_{n}}$ have pairwise disjoint
representatives. The main results of the paper are the following:

\begin{theorem} \label{theorem1} Suppose that $g \geq 2$
and $R$ has at most $g-1$ boundary components. Then a simplicial
map $\lambda : \mathcal{N}(R) \rightarrow \mathcal{N}(R)$ is
superinjective if and only if $\lambda$ is induced by a
homeomorphism of $R$.
\end{theorem}

\begin{theorem} \label{theorem2} Suppose that $g \geq 2$. If $R$ is not a closed
surface of genus two, then $Aut(\mathcal{N}(R))= Mod_R^*$. If $R$
is a closed surface of genus 2, then $Aut(\mathcal{N}(R))= Mod_R
^* /\mathcal{C}(Mod_R ^*)$.
\end{theorem}

\noindent \rule{1.5in}{0.3pt}

\noindent \small{Supported by a Rackham Faculty Fellowship, The
Rackham Graduate School, University of Michigan.}

\begin{theorem} \label{theorem3} Suppose $g = 2$ and $p \leq 1$.
A simplicial map $\lambda : \mathcal{C}(R) \rightarrow
\mathcal{C}(R)$ is superinjective if and only if $\lambda$ is
induced by a homeomorphism of $R$.
\end{theorem}

\begin{theorem}
\label{theorem4} Let $K$ be a finite index subgroup of $Mod_R^*$
and $f$ be an injective homomorphism $f:K \rightarrow Mod_R^*$. If
$g = 2$ and $p=1$ then $f$ has the form $k \rightarrow hkh^{-1}$
for some $h \in Mod_R^*$ and $f$ has a unique extension to an
automorphism of $Mod_R^*$. If $R$ is a closed surface of genus 2,
then $f$ has the form $k \rightarrow hkh^{-1} i^{m(k)}$ for some
$h \in Mod_R^*$ where $m$ is a homomorphism $K \rightarrow
\mathbb{Z}_2$ and $i$ is the hyperelliptic involution on $R$.
\end{theorem}

In section 2, we give some properties of the superinjective
simplicial maps of $\mathcal{N}(R)$. In section 3, we give some
properties of the superinjective simplicial maps of
$\mathcal{C}(R)$ and we prove Theorem \ref{theorem3} and Theorem
\ref{theorem4}. In section 4, we prove that a superinjective
simplicial map $\lambda : \mathcal{N}(R) \rightarrow
\mathcal{N}(R)$ extends to a superinjective simplicial map on
$\mathcal{C}(R)$, and we prove Theorem \ref{theorem1} and Theorem
\ref{theorem2}.\\

Theorem \ref{theorem3} and Theorem \ref{theorem4} complete the
author's previous results given in \cite{Ir1} and \cite{Ir2} to
connected orientable surfaces of genus at least two. These
theorems were motivated by the work of Ivanov \cite{Iv1}, and
Ivanov and McCarthy \cite{IMc}. For automorphism groups of some
complexes related to complex of nonseparating curves
$\mathcal{N}(R)$, see \cite{Luo}, \cite{Sc}.

\section{Properties of Superinjective Simplicial Maps of $\mathcal{N}(R)$}

\noindent A \textit{circle} on $R$ is a properly embedded image of
an embedding $S^{1} \rightarrow R$. A circle on $R$ is said to be
\textit{nontrivial} (or \textit{essential}) if it doesn't bound a
disk and it is not homotopic to a boundary component of $R$. Let
$C$ be a collection of pairwise disjoint circles on $R$. The
surface obtained from $R$ by cutting along $C$ is denoted by
$R_C$. A nontrivial circle $a$ on $R$ is called
\textit{nonseparating} if the surface $R_{a}$ is connected. Let
$\alpha$ and $\beta$ be two vertices in $\mathcal{N}(R)$. The
\textit{geometric intersection number} $i(\alpha, \beta)$ is
defined to be the minimum number of points of $a \cap b$ where $a
\in \alpha$ and $b \in \beta$. A simplicial map $\lambda :
\mathcal{N}(R) \rightarrow \mathcal{N}(R)$ is called {\bf
superinjective} if the following condition holds: if $\alpha,
\beta$ are two vertices in $\mathcal{N}(R)$ such that
$i(\alpha,\beta) \neq 0$, then $i(\lambda(\alpha),\lambda(\beta))
\neq 0$.

\begin{lemma}
\label{injective} Suppose $g \geq 2$ and $p \geq 0$. A
superinjective simplicial map $\lambda : \mathcal{N}(R)
\rightarrow \mathcal{N}(R)$ is injective.
\end{lemma}

\begin{proof} Let $\alpha$ and $\beta$ be two distinct vertices
in $\mathcal{N}(R)$. If $i(\alpha, \beta) \neq 0$, then
$i(\lambda(\alpha), \lambda(\beta)) \neq 0$, since $\lambda$
preserves nondisjointness. So, $\lambda(\alpha) \neq
\lambda(\beta)$. If $i(\alpha, \beta) = 0$, then, since $g \geq 2$
and $p \geq 0$, we can choose a vertex $\gamma$ of
$\mathcal{N}(R)$ such that $i(\gamma, \alpha)= 0$ and $i(\gamma,
\beta) \neq 0$. Then $i(\lambda(\gamma), \lambda(\alpha)) = 0$ and
$i(\lambda(\gamma), \lambda(\beta)) \neq 0$. So, $\lambda(\alpha)
\neq \lambda(\beta)$. Hence $\lambda$ is
injective. \end{proof}\\

Let $P$ be a set of pairwise disjoint circles on $R$. $P$ is
called a {\it pair of pants decomposition} of $R$, if $R_P$ is a
disjoint union of genus zero surfaces with three boundary
components, pairs of pants. A pair of pants of a pants
decomposition is the image of one of these connected components
under the quotient map $q:R_P \rightarrow R$. The image of the
boundary of this component is called the \textit{boundary of the
pair of pants}. A pair of pants is called \textit{embedded} if the
restriction of $q$ to the corresponding component of $R_P$ is an
embedding. An ordered set $(a_1, ..., a_{3g-3+p})$ is called an
{\it ordered pair of pants decomposition} of $R$ if $\{a_1, ...,
a_{3g-3+p}\}$ is a pair of pants decomposition of $R$.

\begin{lemma}
\label{imageofpantsdecomp} Suppose $g \geq 2$ and $p \geq 0$. Let
$\lambda : \mathcal{N}(R) \rightarrow \mathcal{N}(R)$ be a
superinjective simplicial map. Let $P$ be a pair of pants
decomposition consisting of nonseparating circles on $R$. Then
$\lambda$ maps the set of isotopy classes of elements of $P$ to
the set of isotopy classes of elements of a pair of pants
decomposition $P'$ of $R$.
\end{lemma}

\begin{proof} The set of isotopy classes of elements of $P$ forms
a top dimensional simplex, $\bigtriangleup$, in $\mathcal{N}(R)$.
Since $\lambda$ is injective, it maps $\bigtriangleup$ to a top
dimensional simplex $\bigtriangleup'$ in $\mathcal{N}(R)$. A set
of pairwise disjoint representatives of the vertices of
$\bigtriangleup'$ is a pair of pants decomposition $P'$ of
$R$.\end{proof}\\

Let $P$ be a pair of pants decomposition of $R$. Let $a$ and $b$
be two distinct elements in $P$. Then $a$ is called {\it adjacent}
to $b$ w.r.t. $P$ iff there exists a pair of pants in
$P$ which has $a$ and $b$ on its boundary.\\

\noindent {\bf Remark}: Let $P$ be a pair of pants decomposition
of $R$. Let $[P]$ be the set of isotopy classes of elements of
$P$. Let $\alpha, \beta \in [P]$. We say that $\alpha$ is adjacent
to $\beta$ w.r.t. $[P]$ if the representatives of $\alpha$ and
$\beta$ in $P$ are adjacent w.r.t. $P$. By Lemma
\ref{imageofpantsdecomp}, $\lambda$ gives a correspondence on the
isotopy classes of elements of pair of pants decompositions
consisting of nonseparating circles on $R$. $\lambda([P])$ is the
set of isotopy classes of elements of a pair of pants
decomposition which corresponds to $P$, under this correspondence.

\begin{lemma}
\label{adjacent} Suppose $g \geq 2$ and $p \geq 0$. Let $\lambda :
\mathcal{N}(R) \rightarrow \mathcal{N}(R)$ be a superinjective
simplicial map. Let $P$ be a pair of pants decomposition
consisting of nonseparating circles on $R$. Then $\lambda$
preserves the adjacency relation for two circles in $P$, i.e. if
$a, b \in P$, $a$ is adjacent to $b$ w.r.t. $P$ and $[a]=\alpha,
[b]=\beta$, then $\lambda(\alpha)$ is adjacent to $\lambda(\beta)$
w.r.t. $\lambda([P])$.\end{lemma}

\begin{proof} Let $P$ be a pair of pants decomposition consisting of
nonseparating circles on $R$. If $g=2$ and $p \leq 1$, then every
element in $\lambda([P])$ is adjacent to any other element in
$\lambda([P])$, so the lemma is clear. For the other cases; Let
$a, b$ be two adjacent circles in $P$ and $[a]=\alpha$,
$[b]=\beta$. By Lemma \ref{imageofpantsdecomp}, we can choose a
pair of pants decomposition, $P'$, such that $\lambda([P])= [P']$.
Let $P_o$ be a pair of pants of $P$, having $a$ and $b$ on its
boundary. $P_o$ is an embedded pair of pants. There are two
possible cases for $P_o$, depending on whether $a$ and $b$ are the
boundary components of another pair of pants or not. For each of
these cases, we show how to choose a circle $c$ which essentially
intersects $a$ and $b$ and does not intersect any other
circle in $P$ in Figure 1.\\

\begin{figure}[htb]
\begin{center}
\epsfxsize=1.97in \epsfbox{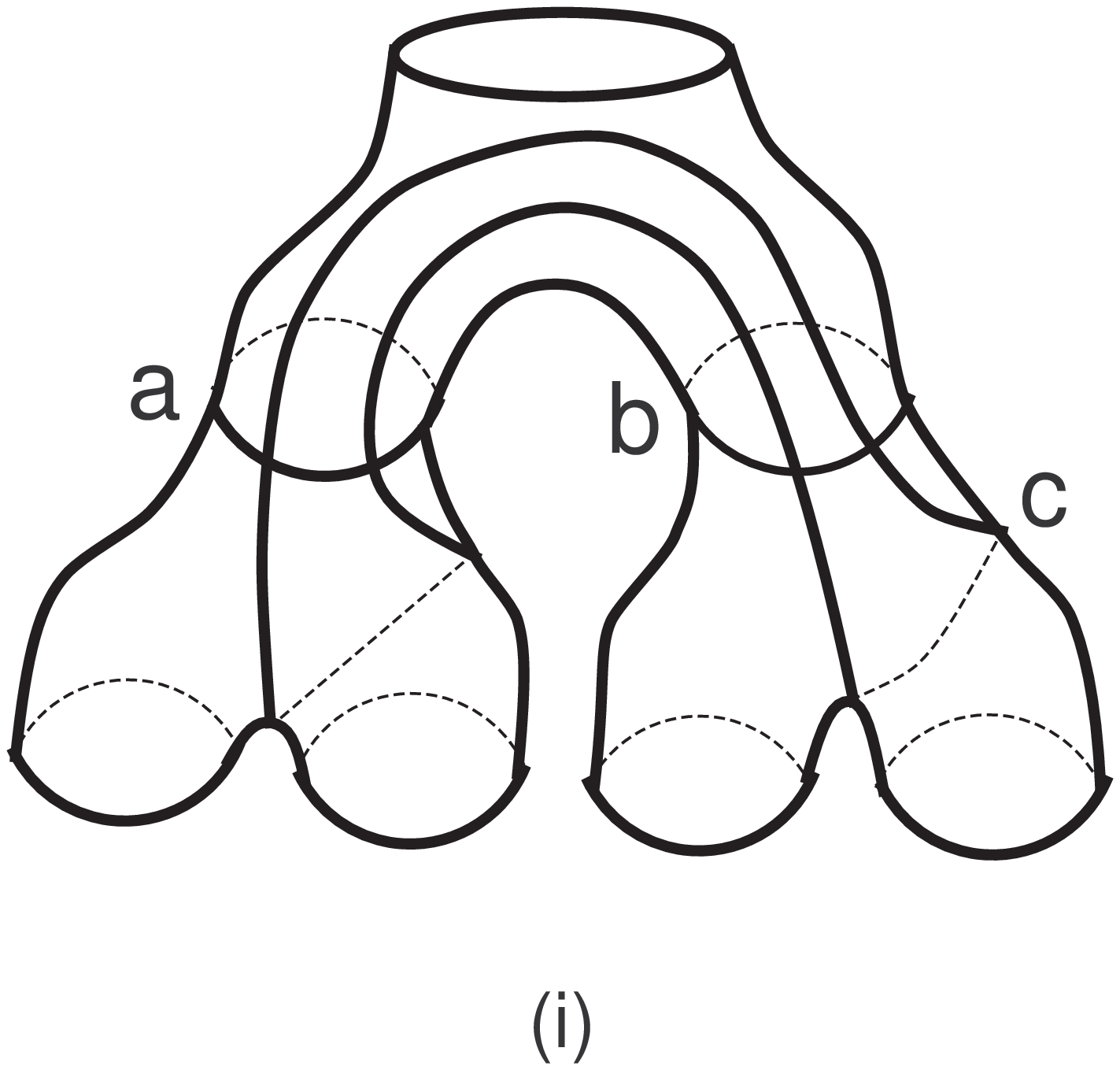} \epsfxsize=1.65in
\epsfbox{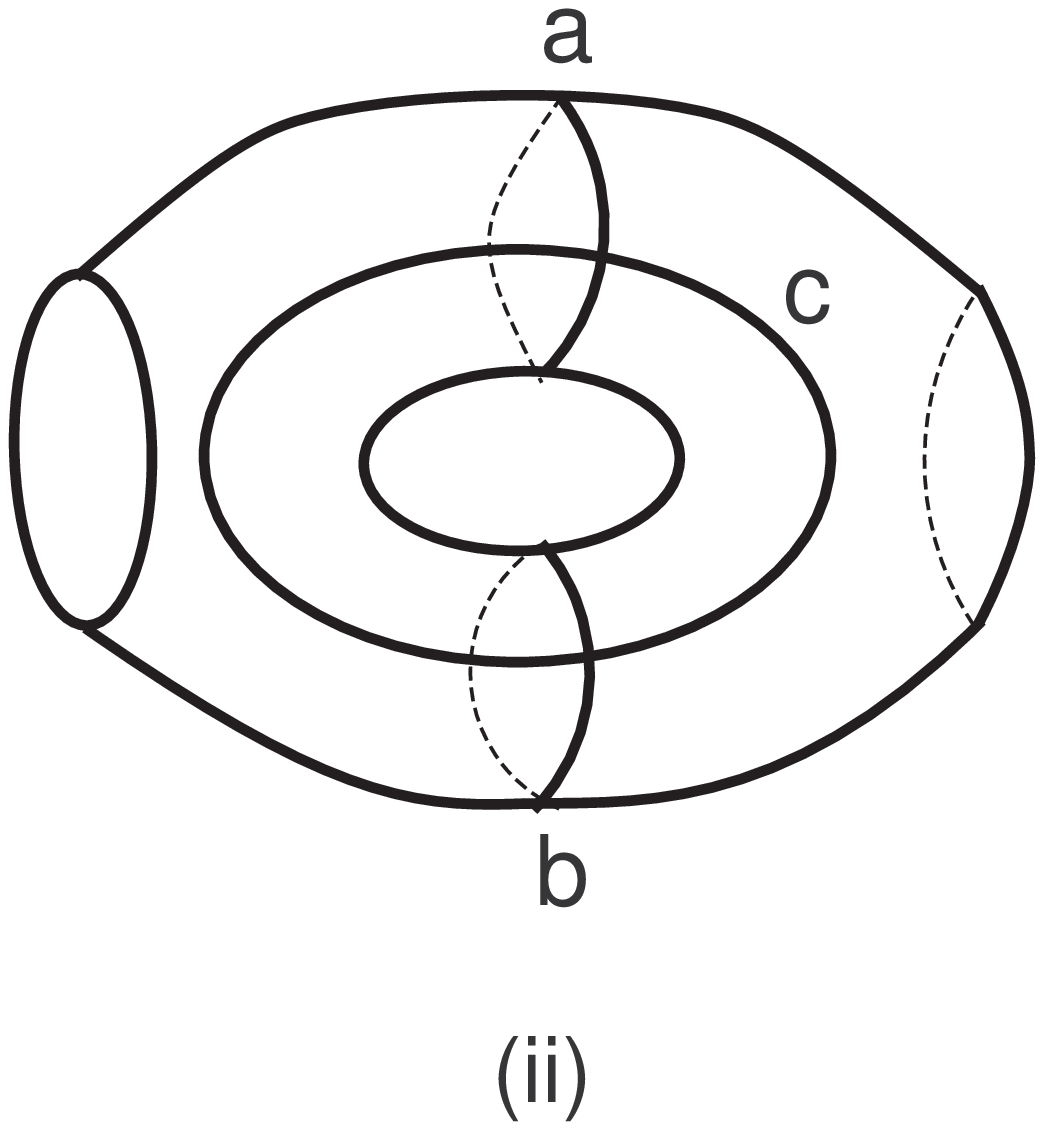} \vspace{0.2cm} \epsfxsize=4.1in
\caption{Two possible cases for $P_o$} \label{picture1}
\end{center}
\end{figure}

Let $\gamma= [c]$. Assume that $\lambda(\alpha)$ and
$\lambda(\beta)$ do not have adjacent representatives. Since
$i(\gamma, \alpha)\neq 0$ and $i(\gamma, \beta) \neq 0$, we have
$i(\lambda (\gamma), \lambda(\alpha)) \neq 0$ and
$i(\lambda(\gamma), \lambda(\beta)) \neq 0$ by superinjectivity.
Since $i(\gamma, [e]) = 0$ for all $e$ in $P \setminus \{a,b\}$,
we have $i(\lambda(\gamma), \lambda([e]))=0$ for all $e$ in $P
\setminus \{a,b\}$. But this is not possible because
$\lambda(\gamma)$ has to intersect geometrically essentially with
some isotopy class other than $\lambda(\alpha)$ and
$\lambda(\beta)$ in the image pair of pants decomposition to be
able to make essential intersections with $\lambda(\alpha)$ and
$\lambda(\beta)$. This is a contradiction to the assumption that
$\lambda(\alpha)$ and $\lambda(\beta)$ do not have adjacent
representatives.\end{proof}\\

Let $P$ be a pair of pants decomposition of $R$. A curve $x \in P$
is called a {\it 4-curve} in $P$ if there exist four distinct
circles in $P$, which are adjacent to $x$ w.r.t. $P$.  Note that
in a pants decomposition every curve is adjacent to at most 4
curves.

\begin{lemma}
\label{embedded} Suppose $g \geq 2$ and $p \geq 0$. Let $\lambda :
\mathcal{N}(R) \rightarrow \mathcal{N}(R)$ be a superinjective
simplicial map and $\alpha, \beta, \gamma$ be distinct vertices in
$\mathcal{N}(R)$ having pairwise disjoint representatives which
bound a pair of pants in $R$. Then $\lambda(\alpha),
\lambda(\beta), \lambda(\gamma)$ are distinct vertices in
$\mathcal{N}(R)$ having pairwise disjoint representatives which
bound a pair of pants in $R$.\end{lemma}

\begin{figure}
\begin{center}
\epsfxsize=5.3in \epsfbox{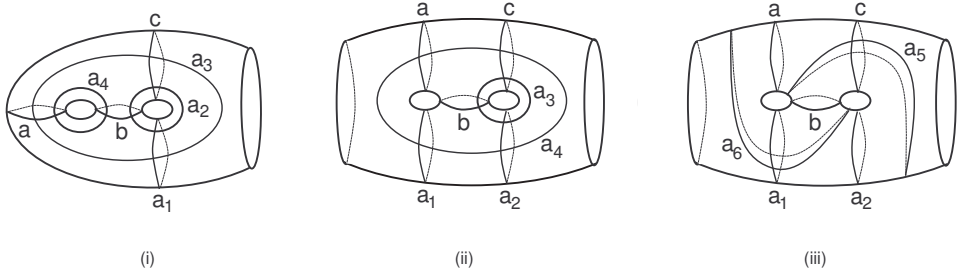} \caption {Curve
configurations}
\end{center}
\end{figure}

\begin{proof} Let $a, b, c$ be pairwise disjoint representatives of
$\alpha, \beta, \gamma$ respectively. If $R$ is a closed surface
of genus two then the statement is obvious. Consider the case when
$g = 2, p = 1$. We complete $\{a, b, c\}$ to a pair of pants
decomposition $P = \{a, b, c, a_1\}$ consisting of nonseparating
circles as shown in Figure 2 (i). Let $P'$ be a pair of pants
decomposition of $R$ such that $\lambda([P]) = [P']$. Let $a', b',
c', a_1'$ be the representatives of $\lambda([a]), \lambda([b]),
\lambda([c]), \lambda([a_1])$ in $P'$ respectively. Since $a$ is
adjacent to $b$ w.r.t. $P$, $a'$ is adjacent to $b'$ w.r.t. $P'$.
Then there is a pair of pants $Q_1$ in $P'$ having $a'$ and $b'$
on its boundary. Let $x$ be the other boundary component of $Q_1$.
If $x = c'$ then we are done. Assume that $x \neq c'$. Then $x$
can't be $a'$, since otherwise $b'$ would be a separating circle
which is a contradiction. Similarly, $x$ can't
be $b'$. Then $x$ is either $a_1'$ or it is the boundary component of $R$.\\

Assume $x$ is a boundary component of $R$. Then, since $a$ is
adjacent to $c, a_1$ w.r.t. $P$, $a'$ is adjacent to $c', a_1'$
w.r.t. $P'$ and so, there is a pair of pants $Q_2$ in $P'$ having
$a'$ and $c', a_1'$ on its boundary. Similarly, since $b$ is
adjacent to $c, a_1$ w.r.t. $P$, $b'$ is adjacent to $c', a_1'$
w.r.t. $P'$ and so, there is a pair of pants $Q_3$ in $P'$ having
$b'$ and $c', a_1'$ on its boundary. Then it is clear that $R =
Q_1 \cup Q_2 \cup Q_3$. Now we consider the circles $a_2, a_3$
which are as shown in Figure 2 (i). Since $a_2$ intersects $b$
essentially and $a_2$ is disjoint from $a$, we can choose a
representative $a_2'$ of $\lambda([a_2])$ such that there exists
an essential arc $w$ of $a_2'$ in $Q_1$ which starts and ends on
$b'$ and which does not intersect $a' \cup \partial R$. Since
$a_3$ is disjoint from $b$ and $a_2$, there exists a
representative $a_3'$ of $\lambda([a_3])$ such that $a_3'$ is
disjoint from $b' \cup w$. But then $a_3'$ could be isotoped so
that it is disjoint from $a'$, since $a'$ is a boundary component
of a regular neighborhood of $b' \cup w$ in $Q_1$. This is a
contradiction since $i([a], [a_3]) \neq 0$ and so $i(\lambda([a]),
\lambda([a_3])) \neq 0$. So, $x$ can not be
a boundary component of $R$.\\

Assume that $x = a_1'$. Then, since $a$ is adjacent to $c$ w.r.t.
$P$, $a'$ is adjacent to $c'$ w.r.t. $P'$ and so, there is a pair
of pants $Q_2$ in $P'$ having $a'$ and $c'$ on its boundary. Let
$y$ be the other boundary component of $Q_2$. If $y = c'$, then
$a'$ would be a separating circle, which is a contradiction. If
$y$ is the boundary component of $R$, then since $c$ is adjacent
to $b, a_1$ w.r.t. $P$, $c'$ is adjacent to $b', a_1'$ in $P'$ and
so, there is a pair of pants $Q_3$ in $P'$ having $b', c', a_1'$
on its boundary. Then it is clear that $R = Q_1 \cup Q_2 \cup
Q_3$. Now we consider the circles $a_2, a_4$ which are as shown in
Figure 2 (i). Since $a_4$ intersects $a$ essentially and $a_4$ is
disjoint from $c$, we can choose a representative $a_4'$ of
$\lambda([a_4])$ such that there exists an essential arc $w$ of
$a_4'$ in $Q_2$ which starts and ends on $a'$ and which does not
intersect $c' \cup \partial R$. Now, we consider $a_2$; $a_2$ is
disjoint from $a$ and $a_4$. Then there exists a representative
$a_2'$ of $\lambda([a_2])$ such that $a_2'$ is disjoint from $a'
\cup w$. But then $a_2'$ could be isotoped so that it is disjoint
from $c'$, since $c'$ is a boundary component of a regular
neighborhood of $a' \cup w$ in $Q_2$. This is a contradiction
since $i([c], [a_2]) \neq 0$ and so $i(\lambda([c]),
\lambda([a_2])) \neq 0$. So, $y$ can not be the boundary component
of $R$. If $y = a_1'$, then there exists a pair of pants $Q_3$
having $b', c'$ and the boundary component of $R$ as its boundary
components. Then we have $R = Q_1 \cup Q_2 \cup Q_3$. Now, we
consider the circles $a_3, a_4$. Since $a_3$ intersects $c$
essentially and $a_3$ is disjoint from $b$, we can choose a
representative $a_3'$ of $\lambda([a_3])$ such that there exists
an essential arc $w$ of $a_3'$ in $Q_3$ which starts and ends on
$c'$ and which does not intersect $b' \cup \partial R$. Now, we
consider $a_4$; $a_4$ is disjoint from $c$ and $a_3$. Then there
exists a representative $a_4'$ of $\lambda([a_4])$ such that
$a_4'$ is disjoint from $c' \cup w$. But then $a_4'$ could be
isotoped so that it is disjoint from $b'$, since $b'$ is a
boundary component of a regular neighborhood of $c' \cup w$ in
$Q_3$. This is contradiction since $i([a_4], [b]) \neq 0$ and so
$i(\lambda([a_4]), \lambda([b])) \neq 0$. Hence, we conclude that
$y = b'$. Then $a', b', c'$ bound a pair of pants
which completes the proof for $g = 2, p = 1$.\\

For the remaining cases we take a subsurface of genus two with two
boundary components, $N$, of $R$, containing nonisotopic circles
$a, b, c, a_1, a_2$ as shown in Figure 2 (ii). Then we complete
$\{a, b, c, a_1, a_2\}$ to a pair of pants decomposition $P$
consisting of nonseparating circles on $R$ in any way we want. Let
$P'$ be a pair of pants decomposition of $R$ such that
$\lambda([P]) = [P']$. Let $a', b', c', a_1', a_2'$ be the
representatives of $\lambda([a]), \lambda([b]), \lambda([c]),
\lambda([a_1]), \lambda([a_2])$ in $P'$ respectively. Since $b$ is
adjacent to $a, c, a_1, a_2$ w.r.t. $P$, $b'$ is adjacent to $a',
c', a_1', a_2'$ w.r.t. $P'$. Then there are two pairs of pants
$Q_1, Q_2$ having $b'$ as boundary components and $Q_1 \cup Q_2$
has $a', c', a_1', a_2'$ on its boundary. W.L.O.G. assume that
$Q_1$ has $a', b'$ on its boundary. Let $x$ be the other boundary
component of $Q_1$. If $x = c'$ then we are done. Assume that $x
\neq c'$. Then $x$ is either $a_1'$ or $a_2'$. Assume that
$x=a_1'$. Since $a_3$ intersects $b$ essentially and $a_3$ is
disjoint from $a \cup a_1$, we can choose a representative $a_3'$
of $\lambda([a_3])$ such that there exists an essential arc $w$ of
$a_3'$ in $Q_1$ which starts and ends on $b'$ and which does not
intersect $a' \cup a_1'$. Now, we consider $a_4$; $a_4$ is
disjoint from $b$ and $a_3$. Then there exists a representative
$a_4'$ of $\lambda([a_4])$ such that $a_4'$ is disjoint from $b'
\cup a_3'$. But then $a_4'$ could be isotoped so that it is
disjoint from $a'$, since $a'$ is a boundary component of a
regular neighborhood of $b' \cup w$ in $Q_1$. This is a
contradiction since $i([a_4], [a]) \neq 0$ and so
$i(\lambda([a_4]), \lambda([a])) \neq 0$. So, $x \neq a_1'$. Now,
assume that $x= a_2'$. We consider the curves $a_5, a_6$ on $N$ as
shown in Figure 2 (iii). Since $a_5$ intersects $c$ essentially
and $a_5$ is disjoint from $a_1, b$, we can choose a
representative $a_5'$ of $\lambda([a_5])$ such that there exists
an essential arc $w$ of $a_5'$ in $Q_2$ which starts and ends on
$c'$ and which does not intersect $a_1' \cup b'$. The circle $a_6$
is disjoint from $c$ and $a_5$. Then there exists a representative
$a_6'$ of $\lambda([a_6])$ such that $a_6'$ is disjoint from $c'
\cup a_5'$. But then $a_6'$ could be isotoped so that it is
disjoint from $a_1'$, since $a_1'$ is a boundary component of a
regular neighborhood of $c' \cup w$ in $Q_2$. This is a
contradiction since $i([a_6], [a_1]) \neq 0$ and so
$i(\lambda([a_6]), \lambda([a_1])) \neq 0$. This completes the
proof of the lemma.\end{proof}\\

Let $\alpha$, $\beta$ be two distinct vertices in
$\mathcal{N}(R)$. We call $(\alpha, \beta)$ to be a
\textit{peripheral pair} in $\mathcal{N}(R)$ if they have disjoint
representatives $x, y$ respectively such that $x, y$ and a
boundary component of $R$ bound a pair of pants in $R$.

\begin{lemma}
\label{peripheral} Suppose $g \geq 2$ and $p \geq 1$. Let $\lambda
: \mathcal{N}(R) \rightarrow \mathcal{N}(R)$ be a superinjective
simplicial map and $(\alpha, \beta)$ be a peripheral pair in
$\mathcal{N}(R)$. Then $(\lambda(\alpha), \lambda(\beta))$ is a
peripheral pair in $\mathcal{N}(R)$.
\end{lemma}

\begin{proof} Let $x, y$ be disjoint representatives of $\alpha, \beta$
respectively such that $x, y$ and a boundary component of $R$
bound a pair of pants in $R$. If $g = 2, p=1$ then we complete $x,
y$ to a pair of pants decomposition $P$ consisting of
nonseparating circles $a, b$ on $R$. Then there are two distinct
pair pants in $P$, one of them has $a, b, x$ on its boundary and
the other has $a, b, y$ on its boundary. Let $P'$ be a pair of
pants decomposition of $R$ such that $\lambda([P]) = [P']$. Let
$x', y', a', b'$ be the representatives of $\lambda([x]),
\lambda([y]), \lambda([a]), \lambda([b]))$ in $P'$ respectively.
By the previous lemma, there exist two pairs of pants $Q_1$, $Q_2$
in $P'$ such that $Q_1$ has $a', b', x'$ on its boundary and $Q_2$
has $a', b', y'$ on its boundary. Then it is clear that $x', y'$
and the boundary component of $R$ bound a pair of pants, which
proves the lemma for this case. If $g \geq 3, p = 1$ then the
proof is similar.\\

If $g = 2, p=2$ then it is easy to see that we can complete $x, y$
to a pair of pants decomposition $P$ consisting of nonseparating
circles $a, b, z$ on $R$ such that $([y], [z])$ is a peripheral
pair, $a, b, x$ bound a pair of pants in $P$ and also $a, b, z$
bound a pair of pants in $P$. Let $P'$ be a pair of pants
decomposition of $R$ such that $\lambda([P]) = [P']$. Let $x', y',
z', a', b'$ be the representatives of $\lambda([x]), \lambda([y]),
\lambda([z]), \lambda([a]), \lambda([b]))$ in $P'$ respectively.
By the previous lemma, there exist two pairs of pants $Q_1$, $Q_2$
in $P'$ such that $Q_1$ has $a', b', x'$ on its boundary and $Q_2$
has $a', b', z'$ on its boundary. Then, since $x$ is adjacent to
$y$ w.r.t. $P$, $x'$ is adjacent to $y'$ w.r.t. $P'$. So, there
exists a pair of pants $Q_3$ in $P'$ having $x', y'$ on its
boundary. Let $w$ be the third boundary component of $Q_3$. It is
easy to see that $Q_1 \cup Q_2 \cup Q_3$ is a genus one surface
with three boundary components, $w, y', z'$. If $w = z'$ then $y'$
would be a separating curve which is a contradiction. If $w = y'$
then $z'$ would be a separating curve which is a contradiction.
Then it is clear that $w$ has to be a boundary component of $R$,
which proves the lemma for this case.\\

Assume that $g = 2, p \geq 3$. We choose distinct essential
circles $z, t, w$ as shown in Figure 3 (i). Then we complete $x,
y, z, t, w$ to a pair of pants decomposition $P$ consisting of
nonseparating circles in any way we like. Let $P'$ be a pair of
pants decomposition of $R$ such that $\lambda([P]) = [P']$. Let
$x', y', z', t', w'$ be the representatives of $\lambda([x]),
\lambda([y]), \lambda([z]), \lambda([t]), \lambda([w])$ in $P'$
respectively. Since $z$ is a 4-curve in $P$, $z'$ is a 4-curve in
$P'$. Since $x, z, w$ are the boundary components of a pair of
pants in $P$, by the previous lemma, $x', z', w'$ are the boundary
components of a pair of pants, $Q_1$, in $P'$. Similarly, since
$z, y, t$ are the boundary components of a pair of pants in $P$,
by the previous lemma, $z', y', t'$ are the boundary components of
a pair of pants, $Q_2$, in $P'$. Since $x$ is adjacent to $y$ in
$P$, $x'$ is adjacent to $y'$ in $P'$. Then there exists a pair of
pants $Q_3$ such that $Q_3$ has $x', y'$ on its boundary. Let $r$
be the boundary component of $Q_3$ different from $x'$ and $y'$.
Suppose that $r$ is not a boundary component of $R$. Then $r$ is
an essential circle. Then $Q_1 \cup Q_2 \cup Q_3$ is a genus one
subsurface with three boundary components $r, w', t'$ which are
nonseparating circles in $R$. Each two of $x', w', t'$ can be
connected by an arc in the complement of $Q_1 \cup Q_2 \cup Q_3$
in $R$. But this would be possible only if genus of $R$ is at
least 3. Since we assumed that $g=2$, we get a contradiction. So,
$r$ has to be a boundary component of $R$, i.e. $(\lambda(\alpha),
\lambda(\beta))$ is a peripheral pair in $\mathcal{N}(R)$.\\

\begin{figure}
\begin{center}
\epsfxsize=3.8in \epsfbox{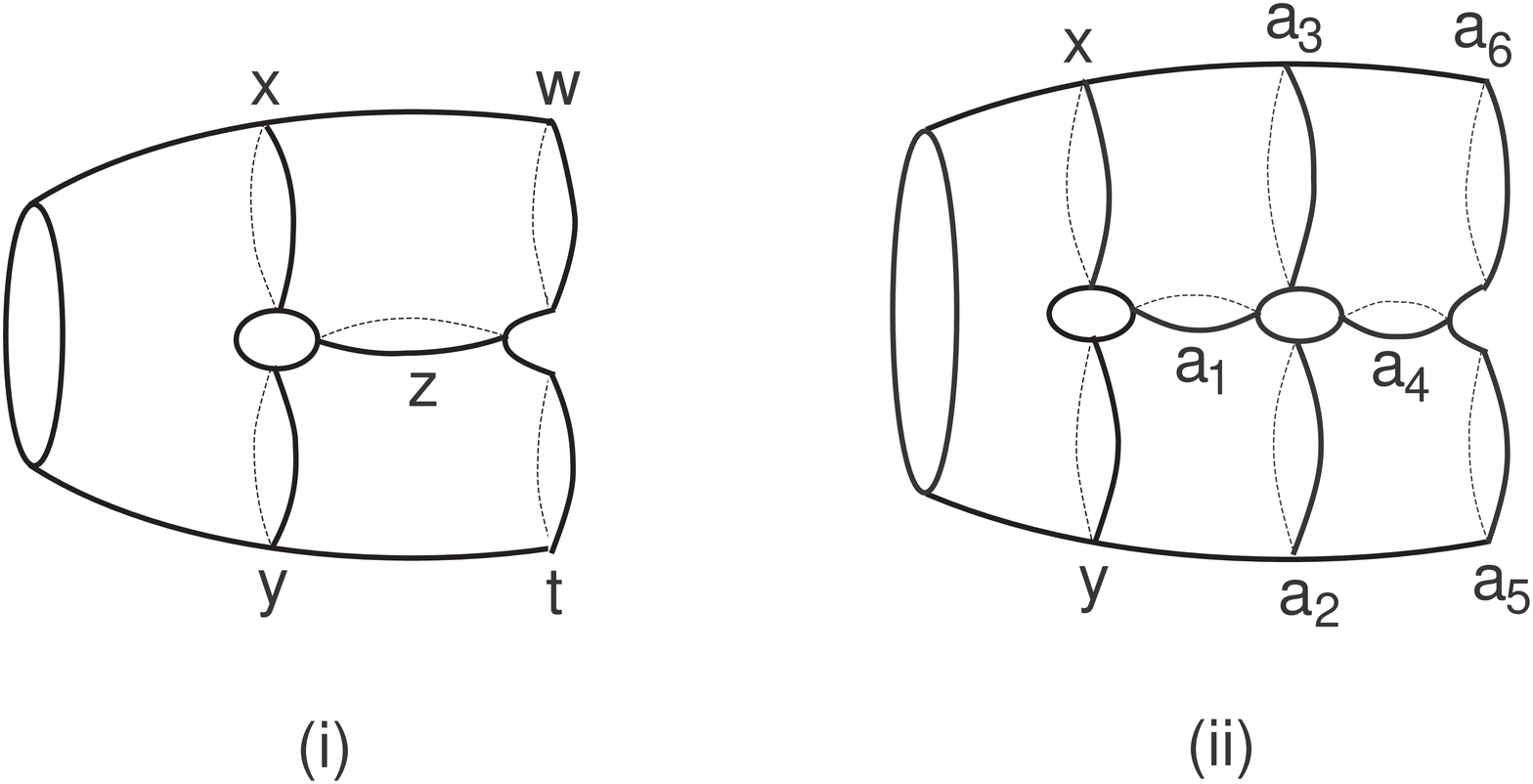} \caption {Pants
decompositions}
\end{center}
\end{figure}

Now, assume that $g=3, p \geq 2$. We choose distinct essential
circles $a_1, ..., a_6$ as shown in Figure 3 (ii). Then we
complete $x, y, a_1, ..., a_6$ to a pair of pants decomposition
$P$ consisting of nonseparating circles in any way we like. Let
$P'$ be a pair of pants decomposition of $R$ such that
$\lambda([P]) = [P']$. Let $x', y', a_1', ..., a_6'$ be the
representatives of $\lambda([x]), \lambda([y]), \lambda([a_1]),
..., \lambda([a_6])$ in $P'$ respectively. Since $a_1$ is a
4-curve in $P$, $a_1'$ is a 4-curve in $P'$. Since $x, a_1, a_3$
are the boundary components of a pair of pants in $P$, by the
previous lemma, $x', a_1', a_3'$ are the boundary components of a
pair of pants, $Q_1$, in $P'$. Similarly, since $y, a_1, a_2$ are
the boundary components of a pair of pants in $P$, by the previous
lemma, $y', a_1', a_2'$ are the boundary components of a pair of
pants, $Q_2$, in $P'$. By using similar arguments we can see that
$a_4'$ is a 4-curve, $a_3', a_4', a_6'$ are the boundary
components of a pair of pants, $Q_3$, in $P'$ and $a_2', a_4',
a_5'$ are the boundary components of a pair of pants, $Q_4$ in
$P'$. Since $x$ is adjacent to $y$ in $P$, $x'$ is adjacent to
$y'$ in $P'$. Then there exists a pair of pants $Q_5$ such that
$Q_5$ has $x', y'$ on its boundary. Let $r$ be the boundary
component of $Q_5$ different from $x'$ and $y'$. Suppose that $r$
is not a boundary component of $R$. Then $r$ is an essential
circle. Then $Q_1 \cup Q_2 \cup Q_3 \cup Q_4 \cup Q_5$ is a genus
two subsurface with three boundary components $r, a_5', a_6'$
which are nonseparating circles in $R$. Each two of $r, a_5',
a_6'$ can be connected by an arc in the complement of $Q_1 \cup
Q_2 \cup Q_3 \cup Q_4 \cup Q_5$ in $R$. But this would be possible
only if genus of $R$ is at least 4. Since we assumed that $g=3$,
we get a contradiction. So, $r$ has to be a boundary component of
$R$, i.e. $(\lambda(\alpha), \lambda(\beta))$ is a peripheral pair
in $\mathcal{N}(R)$. The proof of the case when $g \geq 4, p \geq
2$ is similar.\end{proof}

\begin{lemma}
\label{top} Suppose $g \geq 2$ and $p \geq 0$. Let $\lambda :
\mathcal{N}(R) \rightarrow \mathcal{N}(R)$ be a superinjective
simplicial map. Then $\lambda$ preserves topological equivalence
of ordered pairs of pants decompositions consisting of
nonseparating circles on $R$, (i.e. for a given ordered pair of
pants decomposition $P=(c_1, c_2, ..., c_{3g-3+p})$ of $R$ where
$[c_i] \in \mathcal{N}(R)$, and a corresponding ordered pair of
pants decomposition $P'=(c_1', c_2', ...,$ $ c_{3g-3+p}')$ of $R$,
where $[c_i']= \lambda([c_i])$ $\forall i= 1, 2, ..., 3g-3+p$,
there exists a homeomorphism $H: R \rightarrow R$ such that
$H(c_i)=c_i'$ $\forall i= 1, 2, ..., 3g-3+p$).\end{lemma}

\begin{proof} Let $P = (c_1, c_2, ..., c_{3g-3+p})$ be an ordered pair of pants
decomposition consisting of nonseparating circles on $R$. Let
$c_i' \in \lambda([c_i])$ such that the elements of $\{c_1', c_2',
..., c_{3g-3+p}'\}$ are pairwise disjoint. Then $P'=(c_1', c_2',
..., c_{3g-3+p}')$ is an ordered pair of pants decomposition of
$R$. Let $(B_1, B_2, ..., B_{m})$ be an ordered set of all the
pairs of pants in $P$. By Lemma \ref{embedded} and Lemma
\ref{peripheral}, there is a corresponding, ``image'', ordered
collection of pairs of pants $(B_1', B_2',..., B_{m}')$. A pair of
pants having three essential curves on its boundary corresponds to
a pair of pants having three essential curves on its boundary. A
pair of pants having an inessential boundary component corresponds
to a pair of pants having an inessential boundary component. Then
by the classification of surfaces, there exists an orientation
preserving homeomorphism $h_i : B_i \rightarrow B_i'$, for all $i
=1,..., m$. We can compose each $h_i$ with an orientation
preserving homeomorphism $r_i$ which switches the boundary
components, if necessary, to get $h_i'= r_i \circ h_i$ to agree
with the correspondence given by $\lambda$ on the boundary
components, (i.e. for each essential boundary component $a$ of
$B_i$ for $i=1,...,m$, $\lambda([q(a)])= [q'(h_i'(a))]$ where $q:
R_P \rightarrow R$ and $q': R_{P'} \rightarrow R$ are the natural
quotient maps). Then for two pairs of pants with a common boundary
component, we can glue the homeomorphisms by isotoping the
homeomorphism of the one pair of pants so that it agrees with the
homeomorphism of the other pair of pants on the common boundary
component. By adjusting these homeomorphisms on the boundary
components and gluing them we get a homeomorphism $h : R
\rightarrow R$ such that $h(c_i)=c_i'$ for all $i=1,2, ...,
3g-3+p$.\end{proof}

\begin{lemma} \label{Irmaklemma} Suppose $g \geq 2$ and $p \geq 0$.
Let $\alpha_{1}$ and $\alpha_{2}$ be two vertices in
$\mathcal{N}(R)$. Then $i( \alpha_{1}, \alpha_{2})=1$ if and only
if there exist isotopy classes $\alpha_{3}, \alpha_{4},
\alpha_{5}, \alpha_{6}, \alpha_{7}$ in $\mathcal{N}(R)$ such that

\indent (i) $i(\alpha_{i}, \alpha_{j})=0$ if and only if the
$i^{th}$ and $j^{th}$ circles on Figure 4 are disjoint.

\indent (ii) $\alpha_{1}, \alpha_{3}, \alpha_{5}, \alpha_{6}$ have
pairwise disjoint representatives $a_1, a_3, a_5, a_6$
respectively such that $a_5 \cup a_6$ divides $R$ into two pieces,
one of these is a torus with two holes, $T$, containing some
representatives of the isotopy classes $\alpha_1, \alpha_2$ and
$a_1, a_3, a_5$ bound a pair of pants in $T$ and $a_1, a_3, a_6$
bound a pair of pants in $T$.\end{lemma}

\begin{figure}
\begin{center}
\epsfxsize=2.6in \epsfbox{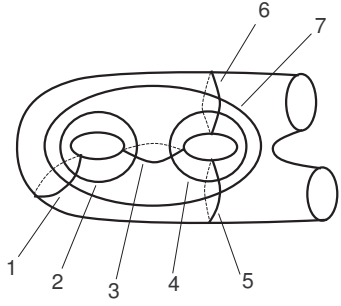} \caption{Circles
intersecting once}
\end{center}
\end{figure}

\begin{proof}Let $i(\alpha_{1}, \alpha_{2})=1$. Let $a_1, a_2$ be
representatives of $\alpha_{1}, \alpha_{2}$ respectively such that
$a_1$ intersects $a_2$ transversely once. Let $N$ be a regular
neighborhood of $a_1 \cup a_2$. Then it is easy to see that $N$ is
a genus one surface with one boundary component and since $R$ is a
surface of genus at least 2, there are simple closed curves as
shown in Figure 4 such that their isotopy classes $\alpha_{3},
\alpha_{4}, \alpha_{5}, \alpha_{6}, \alpha_{7}$ satisfy the
properties (i) and (ii). Suppose that there exist isotopy classes
$\alpha_{1}, \alpha_{2}, ..., \alpha_{7}$ in $\mathcal{N}(R)$
satisfying (i) and (ii). Then we have $a_1, a_3, a_5, a_6$ as
pairwise disjoint representatives of $\alpha_{1}, \alpha_{3},
\alpha_{5}, \alpha_{6}$ respectively such that $a_5 \cup a_6$
divides $R$ into two pieces, one of these is a torus with two
holes, $T$, containing some representatives of the isotopy classes
$\alpha_1, \alpha_2$ and $a_1, a_3, a_5$ bound a pair of pants $P$
in $T$, and $a_1, a_3, a_6$ bound a pair of pants $Q$ in $T$. Let
$a_2, a_4, a_7$ be representatives of $\alpha_{2}, \alpha_{4},
\alpha_{7}$ such that all the curves $a_i$, $i=1,...,7$ have
minimal intersection with each other. Then we have $a_4 \cap a_1 =
\emptyset, a_2 \cap a_6 = \emptyset, a_7 \cap a_3 = \emptyset$.
Since $i(\alpha_{4}, \alpha_{1}) = 0$, $i(\alpha_{4}, \alpha_{3})
\neq 0$ and $i(\alpha_{4}, \alpha_{6}) \neq 0$, $a_4 \cap Q$ has
an arc which connects $a_3$ to $a_6$. Since $i(\alpha_{7},
\alpha_{3}) = 0$, $i(\alpha_{7}, \alpha_{1}) \neq 0$ and
$i(\alpha_{7}, \alpha_{6}) \neq 0$, $a_7 \cap Q$ has an arc which
connects $a_1$ to $a_6$. Similarly, since $i(\alpha_{2},
\alpha_{6}) = 0$, $i(\alpha_{2}, \alpha_{1}) \neq 0$ and
$i(\alpha_{2}, \alpha_{3}) \neq 0$, $a_2 \cap Q$ has an arc which
connects $a_1$ to $a_3$. Then, since $a_2, a_4, a_7$  are pairwise
disjoint, we can see that all the arcs of $a_4 \cap Q$ connect
$a_3$ to $a_6$, all the arcs of $a_7 \cap Q$ connect $a_1$ to
$a_6$ and all the arcs of $a_2 \cap Q$ connect $a_1$ to $a_3$. By
using similar arguments, we can see that all the arcs of $a_4 \cap
P$ connect $a_3$ to $a_5$, all the arcs of $a_7 \cap P$ connect
$a_1$ to $a_5$, and all the arcs of $a_2 \cap P$ connect $a_1$ to
$a_3$. Then by looking at the gluing between the arcs in $a_2 \cap
Q$ and the arcs in $a_2 \cap P$ to form $a_2$, we see that $a_2
\cap Q$ has one arc and $a_2 \cap P$ has one arc. Hence,
$i(\alpha_{1}, \alpha_{2}) = 1$.\end{proof}\\

A characterization of geometric intersection one property in
$\mathcal{C}(R)$ was given by Ivanov, in Lemma 1 in \cite{Iv1}.

\begin{lemma} \label{intone} Suppose $g \geq 2$ and $p \geq
0$. Let $\lambda : \mathcal{N}(R) \rightarrow \mathcal{N}(R)$ be a
superinjective simplicial map. Let $\alpha$, $\beta$ be two
vertices of $\mathcal{N}(R)$. If $i(\alpha, \beta)=1$, then
$i(\lambda(\alpha), \lambda(\beta))=1$.\end{lemma}

\begin{proof} Let $\alpha$, $\beta$ be two vertices of $\mathcal{N}(R)$
such that $i(\alpha, \beta)=1$. Then by Lemma \ref{Irmaklemma},
there exist isotopy classes $\alpha_{3}, \alpha_{4}, \alpha_{5},
\alpha_{6}, \alpha_{7}$ in $\mathcal{N}(R)$ such that
$i(\alpha_{i}, \alpha_{j})=0$ if and only if $i^{th}, j^{th}$
circles on Figure 4 are disjoint and $\alpha_{1}, \alpha_{3},
\alpha_{5}, \alpha_{6}$ have pairwise disjoint representatives
$a_1, a_3, a_5, a_6$ respectively such that $a_5 \cup a_6$ divides
$R$ into two pieces, one of these is a torus with two holes, $T$,
containing some representatives of the isotopy classes $\alpha_1,
\alpha_2$ and $a_1, a_3, a_5$ bound a pair of pants in $T$ and
$a_1, a_3, a_6$ bound a pair of pants in $T$. Then, since
$\lambda$ is superinjective, $i(\lambda(\alpha_{i}),
\lambda(\alpha_{j}))=0$ if and only if $i^{th}, j^{th}$ circles on
Figure 4 are disjoint, and by using Lemma \ref{top} and the
properties that $\lambda$ preserves disjointness and
nondisjointness, we can see that there are pairwise disjoint
representatives $a_1', a_3', a_5', a_6'$ of $\lambda(\alpha_{1}),
\lambda(\alpha_{3}), \lambda(\alpha_{5}), \lambda(\alpha_{6})$
respectively, such that $a_5' \cup a_6'$ divides $R$ into two
pieces, one of these is a torus with two holes, $T$, containing
some representatives of the isotopy classes $\lambda(\alpha_1),
\lambda(\alpha_2)$ and $a_1', a_3', a_5'$ bound a pair of pants in
$T$ and $a_1', a_3', a_6'$ bound a pair of pants in $T$. Then by
Lemma \ref{Irmaklemma}, $i(\lambda(\alpha),
\lambda(\beta))=1$.\end{proof}

\section{Superinjective Simplicial Maps of $\mathcal{C}(R)$ and
Injective Homomorphisms of Finite Index Subgroups of $Mod_R^*$}

If $g = 2$, $p \geq 2$ or $g \geq 3$, $p \geq 0$, by the results
given in \cite{Ir1} and \cite{Ir2}, we have the following: A
simplicial map $\lambda : \mathcal{C}(R) \rightarrow
\mathcal{C}(R)$ is superinjective if and only if $\lambda$ is
induced by a homeomorphism of $R$. If $K$ is a finite index
subgroup of $Mod_R^*$ and $f:K \rightarrow Mod_R^*$ is an
injective homomorphism, then $f$ is induced by a homeomorphism of
$R$. In this section we prove similar results, Theorem
\ref{theorem3} and Theorem \ref{theorem4}, when $g = 2, p \leq 1$.\\

\noindent If $g=2$ and $p \leq 1$, and $\lambda : \mathcal{C}(R)
\rightarrow \mathcal{C}(R)$ is a superinjective simplicial map,
then $\lambda$ is an injective simplicial map which maps pair of
pants decompositions of $R$ to pair of decompositions of $R$ and
it preserves adjacency relation. The proofs of these are similar
to the proofs given in Lemma 1.1-1.3. Now, we prove the following
lemma:

\begin{lemma}
\label{embedded2} Suppose $g=2$ and $p \leq 1$. Let $\lambda :
\mathcal{C}(R) \rightarrow \mathcal{C}(R)$ be a superinjective
simplicial map and $\alpha, \beta, \gamma$ be distinct vertices in
$\mathcal{C}(R)$ having pairwise disjoint representatives which
bound a pair of pants in $R$. Then $\lambda(\alpha),
\lambda(\beta), \lambda(\gamma)$ are distinct vertices in
$\mathcal{C}(R)$ having pairwise disjoint representatives which
bound a pair of pants in $R$.\end{lemma}

\begin{proof} Let $a, b, c$ be pairwise disjoint representatives of
$\alpha, \beta, \gamma$ respectively. If $R$ is a closed surface
of genus two then $a, b, c$ is a pair of pants decomposition on
$R$ and since $a, b, c$ bound a pair of pants, they have to be all
nonseparating circles in this case. Let $P'$ be a pair of pants
decomposition of $R$ such that $\lambda([P]) = [P']$. Let $a', b',
c'$ be the representatives of $\lambda([a]), \lambda([b]),
\lambda([c])$ in $P'$ respectively. Since $a$ is adjacent to $b$
and $c$ w.r.t. $P$, $a'$ is adjacent to $b'$ and $c'$ w.r.t. $P'$.
If one of $a', b'$ or $c'$ was a separating curve on $R$ then the
other two wouldn't be adjacent to each other w.r.t. $P'$, which
would give a contradiction. So, each of $a', b', c'$ is a
nonseparating curve. Then clearly they bound a pair of pants on
$R$. Now assume that $g = 2, p = 1$. There are two cases: either
each of $a, b, c$ is a nonseparating circle or exactly one
of $a, b, c$ is a separating circle.\\


\begin{figure}
\begin{center}
\epsfxsize=5.3in \epsfbox{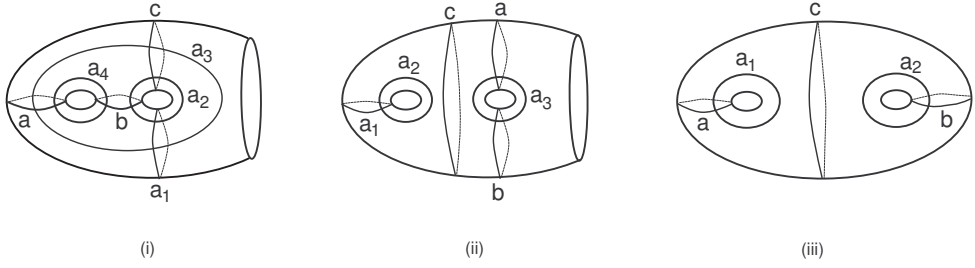} \caption {Curve
configurations}
\end{center}
\end{figure}

Case i: Assume that each of $a, b, c$ is a nonseparating circle
and complete $\{a, b, c\}$ to a pair of pants decomposition $P =
\{a, b, c, a_1\}$ consisting of nonseparating circles as shown in
Figure 5 (i). Let $P'$ be a pair of pants decomposition of $R$
such that $\lambda([P]) = [P']$. Let $a', b', c', a_1'$ be the
representatives of $\lambda([a]), \lambda([b]), \lambda([c]),
\lambda([a_1])$ in $P'$ respectively. Notice that any two curves
in $P$ are adjacent w.r.t. $P$. Since adjacency is preserved, any
two curves in $P'$ must be adjacent w.r.t. $P'$. If one of $a',
b', c', a_1'$ was a separating curve then there would be two
circles in $P'$ which are not adjacent. So, we conclude that all
of $a', b', c', a_1'$ are nonseparating. Then we consider the
curve configuration as shown in Figure 5 (i) and the proof of the
lemma in this case follows as in the proof of Lemma \ref{embedded}
for $g=2, p=1$ case.\\

Case ii: Assume that exactly one of $a, b, c$ is a separating
circle. W.L.O.G. assume that $c$ is a separating circle. Then $c$
separates $R$ into two subsurfaces $R_1, R_2$ as shown in  Figure
5 (ii). Let $a_1, a_2, a_3$ be as shown in the figure. Then $a, b,
c, a_1$ is a pants decomposition $P$ on $R$. Let $P'$ be a pair of
pants decomposition of $R$ such that $\lambda([P]) = [P']$. Let
$a', b', c', a_1'$ be the representatives of $\lambda([a]),
\lambda([b]), \lambda([c]), \lambda([a_1])$ in $P'$ respectively.
Every nonseparating circle $x$ on $R$ could be put inside of a
pair of pants decomposition consisting of nonseparating circles,
and using the method given in case (i), we could see that
$\lambda([x])$ has a nonseparating representative. So, since $a,
b, a_1$ are nonseparating, $a', b', a_1'$ are nonseparating. Since
$a$ is adjacent to $c$ w.r.t. $P$, $a'$ is adjacent to $c'$ w.r.t.
$P'$. Then there is a pair of pants $Q_1$ having $a'$ and $c'$ on
its boundary. Let $x$ be the other boundary component of $Q_1$.
Since $a'$ is not a separating circle, $x$ can't be $c'$. If $x$
is the boundary component of $R$, then since $c'$ is adjacent to
$a_1'$ and $b'$, there exists a pair of pants $Q_2$ in $P'$ having
$c'$, $b'$ and $a_1'$ on its boundary. Then, since $a_2$
intersects $a_1$ essentially and $a_2$ is disjoint from $b$ and
$c$, we can choose a representative $a_2'$ of $\lambda([a_2])$
such that there exists an essential arc $w$ of $a_2'$ in $Q_2$
which starts and ends on $a_1'$ and which does not intersect $b'
\cup c'$. Now, we consider $a_3$ which is disjoint from $a_1$ and
$a_2$. Then there exists a representative $a_3'$ of
$\lambda([a_3])$ such that $a_3'$ is disjoint from $a_1' \cup w$.
But then $a_3'$ could be isotoped so that it is disjoint from
$b'$, since $b'$ is a boundary component of a regular neighborhood
of $a_1' \cup w$ in $Q_2$. This is a contradiction since $i([b],
[a_3]) \neq 0$ and so $i(\lambda([b]), \lambda([a_3])) \neq 0$.
So, $x$ can not be the boundary component of $R$. Then there are
three possibilities: $x$ is one of $a', b', a_1'$. If $x = a'$
then $a'$ wouldn't be adjacent to $b'$ w.r.t. $P'$ which is a
contradiction. Assume $x = a_1'$. Then, since $a_2$ intersects
$a_1$ essentially and $a_2$ is disjoint from $a$ and $c$, we can
choose a representative $a_2'$ of $\lambda([a_2])$ such that there
exists an essential arc $w$ of $a_2'$ in $Q_1$ which starts and
ends on $a_1'$ and which does not intersect $a' \cup c'$. Now, we
consider $a_3$ which is disjoint from $a_1$ and $a_2$. Then there
exists a representative $a_3'$ of $\lambda([a_3])$ such that
$a_3'$ is disjoint from $a_1' \cup w$. But then $a_3'$ could be
isotoped so that it is disjoint from $a'$, since $a'$ is a
boundary component of a regular neighborhood of $a_1' \cup w$ in
$Q_1$. This is a contradiction since $i([a ], [a_3]) \neq 0$ and
so $i(\lambda([a]), \lambda([a_3])) \neq 0$. So, $x \neq a_1'$.
Then $x= b'$, and $a', b', c'$ bound a pair of pants in $P'$.\end{proof}\\

Let $\alpha$, $\beta$ be two distinct vertices in
$\mathcal{C}(R)$. We call $(\alpha, \beta)$ to be a
\textit{peripheral pair} in $\mathcal{C}(R)$ if they have disjoint
representatives $x, y$ respectively such that $x, y$ and a
boundary component of $R$ bound a pair of pants in $R$.

\begin{lemma}
\label{peripheral2} Suppose $g = 2$ and $p=1$. Let $\lambda :
\mathcal{C}(R) \rightarrow \mathcal{C}(R)$ be a superinjective
simplicial map and $(\alpha, \beta)$ be a peripheral pair in
$\mathcal{C}(R)$. Then $(\lambda(\alpha), \lambda(\beta))$ is a
peripheral pair in $\mathcal{C}(R)$.
\end{lemma}

\begin{proof} Let $x, y$ be disjoint representatives of $\alpha, \beta$
respectively such that $x, y$ and a boundary component of $R$
bound a pair of pants in $R$. There are two cases to consider:
Case i: $x$ and $y$ are both nonseparating. Case ii: $x$ and $y$
are both separating. The proof in the first case is similar to the
proof given in Lemma \ref{peripheral}. For the second case, we
complete $x, y$ to a pair of pants decomposition $Q$ consisting of
nonseparating circles $a, b$ on $R$ such that $a$ is in the torus
with one hole which comes with the separation by $x$. Then we will
replace $y$ with a nonseparating curve $w$ such that $a, x, b, w$
is a pants decomposition $P$ on $R$, $x, w, b$ bound a pair of
pants and $w, b$ and the boundary component of $R$ bound a pair of
pants on $R$. Let $P'$ be a pair of pants decomposition of $R$
such that $\lambda([P]) = [P']$. Let $x', w', a', b'$ be the
representatives of $\lambda([x]), \lambda([w]), \lambda([a]),
\lambda([b]))$ in $P'$ respectively. Then by using case i and
Lemma \ref{embedded2}, we see that $x'$ is a separating circle of
genus one. Similarly, $y'$ is a separating circle of genus one on
$R$. Then it is easy to see that $x', y'$ and the boundary
component of $R$ bound a pair of pants.\end{proof}

\begin{lemma}
\label{top2} Suppose $g = 2$ and $p \leq 1$. Let $\lambda :
\mathcal{C}(R) \rightarrow \mathcal{C}(R)$ be a superinjective
simplicial map. Then $\lambda$ preserves topological equivalence
of ordered pairs of pants decompositions on $R$, (i.e. for a given
ordered pair of pants decomposition $P=(c_1, c_2, ..., c_{3+p})$
of $R$ where $[c_i] \in \mathcal{C}(R)$, and a corresponding
ordered pair of pants decomposition $P'=(c_1', c_2', ...,$ $
c_{3+p}')$ of $R$, where $[c_i']= \lambda([c_i])$ $\forall i= 1,
2, ..., 3 + p$, there exists a homeomorphism $H: R \rightarrow R$
such that $H(c_i)=c_i'$ $\forall i= 1, 2, ..., 3+p$).\end{lemma}

\begin{proof} First we will consider the case when $g = 2$ and $p=0$.
Let $P$ be a pair of pants decomposition of $R$ and $A$ be a
nonembedded pair of pants in $P$. The boundary of $A$ consists of
the circles $a, c$ where $c$ is a 1-separating circle on $R$ and
$a$ is a nonseparating circle on $R$. Let $b, a_1, a_2$ be as
shown in Figure 5 (iii). Then $P= \{a, b, c\}$ is a pants
decomposition on $R$. Let $P'$ be a pair of pants decomposition of
$R$ such that $\lambda([P]) = [P']$. Let $a', b', c'$ be the
representatives of $\lambda([a]), \lambda([b]), \lambda([c])$ in
$P'$ respectively. Every nonseparating circle $x$ on $R$ could be
put inside of a pair of pants decomposition consisting of
nonseparating circles, and using the method given in case (i) in
Lemma \ref{embedded2}, we could see that $\lambda([x])$ has a
nonseparating representative. So, since $a, b$ are nonseparating,
$a', b'$ are nonseparating. Assume that $c'$ is also
nonseparating. Then there exist two pairs of pants $Q_1, Q_2$ in
$P'$ such that both of them have $a', b', c'$ on their boundary.
Then, since $a_1$ intersects $a$ essentially and $a_1$ is disjoint
from $b$ and $c$, we can choose a representative $a_1'$ of
$\lambda([a_1])$ such that there exists an essential arc $w$ of
$a_1'$ in $Q_1$ which starts and ends on $a'$ and which does not
intersect $b' \cup c'$. Now, we consider $a_2$ which is disjoint
from $a$ and $a_1$. There exists a representative $a_2'$ of
$\lambda([a_2])$ such that $a_2'$ is disjoint from $a_1' \cup a'$.
But then $a_2'$ could be isotoped so that it is disjoint from
$b'$, since $b'$ is a boundary component of a regular neighborhood
of $a' \cup w$ in $Q_1$. This is a contradiction since $i([b],
[a_2]) \neq 0$ and so $i(\lambda([b]), \lambda([a_2])) \neq 0$.
So, $c'$ has to be separating. Then clearly $a'$ and $c'$ are the
boundary components of a nonembedded pair of pants. Hence,
nonembedded pair of pants in $P$ corresponds to a nonembedded pair
of pants in $P'$. When $g = 2$ and $p=1$, this result follows from
Lemma \ref{embedded2} and Lemma \ref{peripheral2} considering the
circles in Figure 5 (ii).\\

Suppose $g = 2$ and $p \leq 1$. Let $P = (c_1, c_2, ..., c_{3+p})$
be an ordered pair of pants decomposition on $R$. Let $c_i' \in
\lambda([c_i])$ such that the elements of $\{c_1', c_2', ...,
c_{3+p}'\}$ are pairwise disjoint. Then $P'=(c_1', c_2', ...,
c_{3+p}')$ is an ordered pair of pants decomposition of $R$. Let
$(B_1, B_2, ..., B_{m})$ be an ordered set of all the pairs of
pants in $P$. By Lemma \ref{embedded2}, Lemma \ref{peripheral2}
and the arguments given above, there is a corresponding,
``image'', ordered collection of pairs of pants $(B_1', B_2',...,
B_{m}')$. Nonembedded pairs of pants correspond to nonembedded
pairs of pants, embedded pairs of pants correspond to embedded
pairs of pants, and a pair of pants having an inessential boundary
component corresponds to a pair of pants having an inessential
boundary component. Then the proof follows as in Lemma \ref{top}.
\end{proof}


\begin{lemma} \label{intone2} Suppose $g = 2$ and $p \leq 1$. Let $\lambda :
\mathcal{C}(R) \rightarrow \mathcal{C}(R)$ be a superinjective
simplicial map. Let $\alpha$, $\beta$ be two vertices of
$\mathcal{C}(R)$. If $i(\alpha, \beta)=1$, then
$i(\lambda(\alpha), \lambda(\beta))=1$.\end{lemma}

\begin{proof} The proof follows as in the proof of Lemma \ref{intone},
using Lemma \ref{top2}.\end{proof}

\begin{lemma}
\label{horver2} Suppose $g = 2$ and $p \leq 1$. Let $\lambda :
\mathcal{C}(R) \rightarrow \mathcal{C}(R)$ be a superinjective
simplicial map. Let $\alpha$ and $\beta$ be two vertices in
$\mathcal{C}(R)$ which have representatives with geometric
intersection 2 and algebraic intersection 0 on $R$. Then
$\lambda(\alpha)$ and $\lambda(\beta)$ have representatives with
geometric intersection 2 and algebraic intersection 0 on $R$.
\end{lemma}

\begin{figure}
\begin{center} \epsfxsize=6in \epsfbox{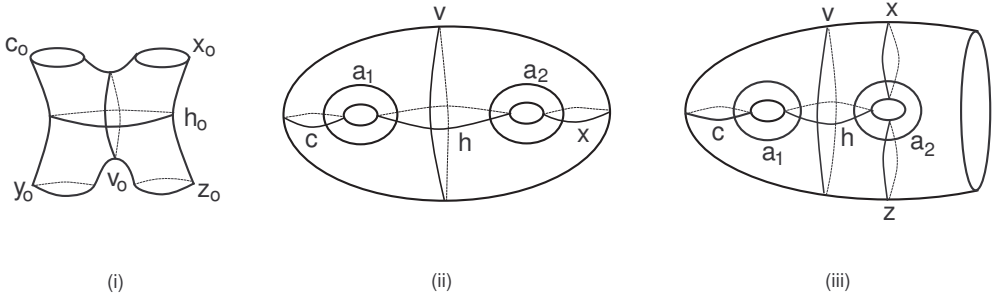} \caption {Curves
intersecting twice}
\end{center}
\end{figure}

\begin{proof} Assume that $g=2$ and $p=0$. Let $h, v$ be representatives
of $\alpha, \beta$ with geometric intersection 2 and algebraic
intersection 0 on $R$. Let $N$ be a regular neighborhood of $h
\cup v$ in $R$. Then $N$ is a sphere with four boundary
components. Let $c, x, y, z$ be boundary components of $N$ such
that there exists a homeomorphism $\varphi : (N, c, x, y, z, h,
v)$ $\rightarrow (N_o, c_o, x_o, y_o, z_o, h_o, v_o)$ where $N_o$
is a standard sphere with four holes having $c_o, x_o, y_o, z_o$
on its boundary and $h_o, v_o$ (horizontal, vertical) are two
circles as indicated in Figure 6 (i). Since $h$ and $v$ have
geometric intersection 2 and algebraic intersection 0 on $R$, none
of $c, x, y, z$ bound a disk on $R$. There are two cases to
consider: either exactly one of $h$ or $v$ is separating or both
$h$ and $v$ are nonseparating. Case i: W.L.O.G. Assume that $v$ is
separating and $h$ is nonseparating. Then $x$ is isotopic to $z$
and $y$ is isotopic to $c$ and $c, x, h, v$ are as shown in Figure
6 (ii). Let $a_1, a_2$ be as shown in the figure. Let $c', h', x'$
be pairwise disjoint representatives of $\lambda([c]),
\lambda([h]), \lambda([x])$ respectively. Then $\{c', h', x'\}$ is
a pair of pants decomposition on $R$. Since nondisjointness and
intersection one is preserved by $\lambda$, there are disjoint
representatives $a_1', a_2'$ of $\lambda([a_1]), \lambda([a_2])$
respectively such that $a_1'$ intersects each of $c'$ and $h'$
exactly once and $a_1'$ is disjoint from $x'$, and similarly
$a_2'$ intersects each of $x'$ and $h'$ exactly once and $a_2'$ is
disjoint from $c'$. Then, since $v$ is disjoint from $c \cup a_1
\cup a_2 \cup x$, we can choose a representative $v'$ of
$\lambda([v])$ such that it lives inside of the cylinder that we
get after cutting $R$ along $c' \cup a_1' \cup a_2' \cup x'$. Then
it is easy to see that $h'$ and $v'$ intersects geometrically
twice and algebraically zero times.\\

\begin{figure}
\begin{center}
\epsfxsize=5.8in \epsfbox{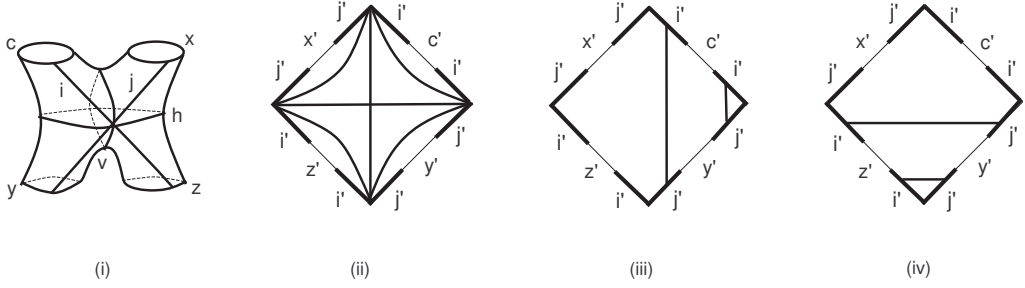} \caption {Curves
intersecting twice}
\end{center}
\end{figure}

Case ii: Assume that both of $h, v$ are nonseparating. Then $c$
and $z$ are the boundary components of an annulus, $A_1$, in $R
\setminus N$ and x and y are the boundary components of an
annulus, $A_2$, in $R \setminus N$. Let $i$ and $j$ be as shown in
Figure 7 (i). We connect the ends points of $i$ on $N$ with an arc
in $A_1$ to get a circle $w$ such that $w$ intersects each of $c$
and $z$ at only one point. Similarly, we connect the ends points
of $j$ on $N$ with an arc in $A_2$ to get a circle $k$ such that
$k$ intersects each of $x$ and $y$ at only one point. Notice that
$w$ and $k$ intersect at only one point. Let $c', x', h', y', z'$
be pairwise disjoint representatives of the image curves. Let $N'$
be sphere with four holes bounded by $c', x', y', z'$. By Lemma
\ref{top2}, $h'$ gives a pants decompositions on $N'$. Now we
choose minimally intersecting representatives $w', k'$ of
$\lambda([w]), \lambda([k])$ respectively such that each of $w'$
and $k'$ intersects $c', x', y', z'$ minimally. Then, since the
intersection one is preserved, there are arcs $i'$, $j'$ of $w'$,
$k'$ in $N'$ respectively such that $i'$ intersects $j'$ once in
$N'$, $i'$ connects $c'$ to $z'$ and $j'$ connects $x'$ to $y'$.
Since $v$ does not intersect any of $x, c, y, z$, and $v$
intersects $w$ and $k$ exactly once, there is a representative
$v'$ of $\lambda([v])$ in $N'$ such that $v'$ intersects each of
$i'$ and $j'$ exactly once. Notice that $h'$ also intersects each
of $i'$ and $j'$ exactly once. When we cut $N'$ along the arcs of
$i'$ and $j'$ we get a disk. The boundary of the disk either is as
shown in Figure 7 (ii) or it has the similar form where only $x'$
and $y'$ are switched. We will consider the first case (the
arguments follow in the second case similarly). If $v'$ makes its
intersection with $i'$ and $j'$ at the intersection point of $i'$
and $j'$ with each other, then $v'$ has to be one of the arcs
shown in Figure 7 (ii), and by looking at the intersection of $v'$
with the other curves, it is easy to see that $v'$ intersects $h'$
geometrically twice and algebraically 0 times on $R$. Suppose that
$v'$ intersects $i'$ and $j'$ at different points. Then there are
two arcs of $v'$ connecting $i'$ to $j'$ in the disk. Then, since
$v'$ does not intersect any of $x', c', y', z'$ we see that $v'$
has to be one of the curves as shown in Figure 7 (iii), (iv).
Since $v'$ is not isotopic to $h'$, $v'$ has to be the curve shown
in Figure 7 (iii) and hence $v'$ intersects $h'$ geometrically
twice and algebraically 0 times on $R$. If $g = 2$ and $p=1$ the
proof is similar (for Case (i) we use Figure 6
(iii)).\end{proof}\\

The proof of Theorem \ref{theorem3} follows from Lemma \ref{top2},
Lemma \ref{intone2}, Lemma \ref{horver2}, and the techniques given
for the proof of Theorem 1.1 in \cite{Ir1}, \cite{Ir2}.\\

A mapping class $g \in Mod_R^*$ is called \textit{pseudo-Anosov}
if $\mathcal{A}$ is nonempty and if $g ^n (\alpha) \neq \alpha$,
for all $\alpha$ in $\mathcal{A}$ and any $n \neq 0$. $g$ is
called \textit{reducible} if there is a nonempty subset $
\mathcal{B} \subseteq \mathcal{A}$ such that a set of disjoint
representatives can be chosen for $\mathcal{B}$ and
$g(\mathcal{B}) = \mathcal{B}$. In this case, $ \mathcal{B}$ is
called a \textit{reduction system} for $g$. Each element of
$\mathcal{B}$ is called a \textit{reduction class} for $g$. A
reduction class, $\alpha$, for $g$, is called an \textit{essential
reduction class} for $g$, if for each $\beta \in \mathcal{A}$ such
that $i(\alpha, \beta) \neq 0$ and for each integer $m \neq 0$,
$g^m (\beta) \neq \beta$. The set, $\mathcal{B}_g$, of all
essential reduction classes for $g$ is called the
\textit{canonical reduction system} for $g$.\\

Let $\Gamma'=ker(\varphi)$ where $\varphi: Mod_R^* \rightarrow
Aut(H_1(R, \mathbb{Z}_3))$ is the homomorphism defined by the
action of homeomorphisms on the homology. The proofs of the Lemma
\ref{rank=1} and Lemma \ref{reducible} follow by the techniques
given in \cite{Ir1}. Note that we need to use that the maximal
rank of an abelian subgroup of $Mod_R^*$ is $3g-3+p$, \cite{BLM}.

\begin{lemma}
\label{rank=1} Suppose $g=2$ and $p \leq 1$. Let $K$ be a finite
index subgroup of $Mod_R^*$ and $f:K \rightarrow Mod_R^*$ be an
injective homomorphism. Let $\alpha \in \mathcal{A}$. Then there
exists $N \in \mathbb{Z^*}$ such that $rank$ $C(C_{\Gamma'}
(f(t_{\alpha} ^{N})) ) = 1$.\end{lemma}

\begin{lemma}
\label{reducible} Suppose $g=2$ and $p \leq 1$. Let K be a finite
index subgroup of $Mod_R^*$. Let $f:K \rightarrow Mod_R^*$ be an
injective homomorphism. Then there exists $N \in \mathbb{Z^*}$
such that $f(t_{\alpha} ^ N)$ is a reducible element of infinite
order for all $\alpha \in \mathcal{A}$.
\end{lemma}

In the proof of Lemma \ref{reducible}, we use that the centralizer
of a p-Anosov element in the extended mapping class group is a
virtually infinite cyclic group, \cite{Mc2}.

\begin{lemma}
\label{correspondence} Suppose $g=2$ and $p \leq 1$. Let $K$ be a
finite index subgroup of $Mod_R^*$ and $f:K \rightarrow Mod_R^*$
be an injective homomorphism. Then $\forall \alpha \in
\mathcal{A}$, $f( t_\alpha ^N)= t_{\beta(\alpha)}^M$ for some $M,
N \in \mathbb{Z^*}$, $\beta(\alpha) \in \mathcal{A}$.
\end{lemma}

\begin{proof} Let $\Gamma= f^{-1}(\Gamma') \cap \Gamma'$. Since
$\Gamma$ is a finite index subgroup we can choose $N \in Z^*$ such
that $t_\alpha^N \in \Gamma$ for all $\alpha$ in $\mathcal{A}$. By
Lemma \ref{reducible} $f(t_{\alpha} ^ N)$ is a reducible element
of infinite order in $Mod_R^*$. Let $C$ be a realization of the
canonical reduction system of $f(t_{\alpha}^N)$. Let $c$ be the
number of components of $C$ and $r$ be the number of p-Anosov
components of $f(t_{\alpha} ^N)$. Since $t_{\alpha} ^ N \in
\Gamma, f(t_{\alpha} ^ N) \in \Gamma'$. By Theorem 5.9 \cite{IMc},
$C(C_{\Gamma'} (f(t_{\alpha} ^ N )))$ is a free abelian group of
rank $c+r$. By Lemma \ref{rank=1} $c+r=1$. Then either $c=1$,
$r=0$ or $c=0$, $r=1$. Since there is at least one curve in the
canonical reduction system we have $c=1$, $r=0$. Hence, since
$f(t_{\alpha} ^ N) \in \Gamma'$, $f(t_{\alpha} ^{N}) = t_{\beta
({\alpha})}^{M}$ for some $M \in \mathbb{Z^*}$, $\beta(\alpha) \in
\mathcal{A}$, \cite{BLM}, \cite{IMc}.\end{proof}\\

\noindent {\bf Remark:} Suppose that $f(t_{\alpha} ^{M}) = t_\beta
^P$ for some $\beta \in \mathcal{A}$ and $M, P \in \mathbb{Z^*}$
and $f(t_{\alpha} ^{N}) = t_\gamma ^Q$ for some $\gamma \in
\mathcal{A}$ and $N, Q \in \mathbb{Z^*}$. Since $f(t_{\alpha} ^{M
\cdot N}) = f(t_{\alpha} ^{N \cdot M})$, $t_\beta ^{PN} = t_\gamma
^{QM}$, $P, Q, M, N \in \mathbb{Z^*}$. Then $\beta = \gamma$.
Therefore, by Lemma \ref{correspondence}, $f$ gives a
correspondence between isotopy classes of circles and $f$ induces
a map, $f_*: \mathcal{A} \rightarrow \mathcal{A}$,
where $f_*(\alpha) = \beta(\alpha)$.\\

In the following lemma we use a well known fact that $f t_\alpha
f^{-1}=t_{f(\alpha)} ^{\epsilon(f)}$ for all $\alpha$ in
$\mathcal{A}$, $f \in Mod_R^*$, where $\epsilon(f) = 1$ if $f$ has
an orientation preserving representative and $\epsilon(f) = -1$ if
$f$ has an orientation reversing representative.

\begin{lemma}
\label{identity} Let $K$ be a finite index subgroup of $Mod_R^*$.
Let $f:K \rightarrow Mod_R^*$ be an injective homomorphism. Assume
that there exists $N \in \mathbb{Z}^*$ such that $\forall \alpha
\in$ $\mathcal{A}$, $\exists Q \in \mathbb{Z}^*$ such that
$f(t_{\alpha} ^N) = t_{\alpha}^Q$. If $g=2$ and $p=1$ then $f$ is
the identity on $K$. If $g=2$ and $p=0$ then $f(k) = ki^{m(k)}$
where $i$ is the hyperelliptic involution on $R$ and $m(k) \in
\{0,1\}$.\end{lemma}

\begin{proof}
We use Ivanov's trick to see that $f(kt_{\alpha} ^ N k^{-1})=$
$f(t_{k(\alpha)} ^{\epsilon{(k)} \cdot N }) = t_{k(\alpha)} ^{Q
\cdot \epsilon{(k)}}$ and $f(kt_{\alpha} ^ N k^{-1}) = f(k)
f(t_{\alpha} ^N) f(k)^{-1}=$ $f(k) t_{\alpha} ^Q f(k)^{-1} =
t_{f(k)(\alpha)} ^{\epsilon(f(k))\cdot Q}$ $\forall \alpha \in
\mathcal{A}$, $\forall k \in K$. Then we have $t_{k(\alpha)} ^{Q
\cdot \epsilon{(k)}} = t_{f(k)(\alpha)} ^{\epsilon(f(k)) \cdot Q}$
$\forall \alpha \in \mathcal{A}$, $\forall k \in K$. Hence,
$k(\alpha) = f(k)(\alpha)$ $\forall \alpha \in \mathcal{A}$,
$\forall k \in K$. Then $k^{-1}f(k)(\alpha) = \alpha$ $\forall
\alpha \in \mathcal{A}$, $\forall k \in K$. Therefore,
$k^{-1}f(k)$ commutes with $t_{\alpha}$ $\forall \alpha \in
\mathcal{A}$, $\forall k \in K$. Then, since $Mod_R$ is generated
by Dehn twists when $g=2, p \leq 1$, $k^{-1}f(k) \in C(Mod_R)$
$\forall k \in K$. If $g=2$ and $p=1$, then $C(Mod_R)$ is trivial
by 5.3 in \cite{IMc}. So, $k = f(k)$ $\forall k \in K$. Hence,
$f=id_K$. If $g=2$ and $p=0$, then $C(Mod_R) = \{id_R, i\} =
\mathbb{Z}_2$, where $i$ is the hyperelliptic involution on $R$.
Then for each $k \in K$ either $k^{-1}f(k)=id_R$ or $k^{-1}f(k) =
i$. So, $f(k) = ki^{m(k)}$ where $m(k) \in \{0,1\}$.\end{proof}

\begin{coroll}
\label{id} Suppose $g=2$ and $p=1$. Let $h: Mod_R^* \rightarrow
Mod_R^*$ be an isomorphism and $f : Mod_R^* \rightarrow Mod_R^*$
be an injective homomorphism.  Assume that there exists $N \in
\mathbb{Z}^*$ such that $\forall \alpha \in$ $\mathcal{A}$,
$\exists Q \in \mathbb{Z}^*$ such that $h(t_{\alpha} ^N) =
f(t_{\alpha}^Q)$. Then $h=f$.
\end{coroll}

\begin{proof} Apply Lemma \ref{identity} to $h^{-1} f$ with $K = Mod_R^*$.
Since for all $\alpha$ in $\mathcal{A}$, $h^{-1} f(t_{\alpha} ^N)
= t_{\alpha} ^{Q}$, we have $h^{-1} f = id_K$. Hence, $h =
f$.\end{proof}\\

By the remark after Lemma \ref{correspondence}, we have that $f: K
\rightarrow Mod_R^*$ induces a map $f_*: \mathcal{A} \rightarrow
\mathcal{A}$, where $K$ is a finite index subgroup of $Mod_R^*$.
In the following lemma we prove that $f_*$ is a superinjective
simplicial map on $\mathcal{C}(R)$.

\begin{lemma}
\label{intersection0} Suppose $g=2$ and $p \leq 1$. Let $f:K
\rightarrow Mod_R^*$ be an injection. Let $\alpha$, $\beta \in
\mathcal{A}$. Then $i(\alpha,\beta)=0 \Leftrightarrow
i(f_{*}(\alpha), f_{*}(\beta))=0$.
\end{lemma}

\begin{proof} There exists $N \in \mathbb {Z^*}$ such that
$t_{\alpha} ^N \in K$ and $t_{\beta} ^N \in K$. Then we have the
following: $i(\alpha, \beta)=0$ $\Leftrightarrow$ $t_{\alpha} ^N
t_\beta ^N = t_{\beta} ^N t_\alpha ^N$ $\Leftrightarrow$
$f(t_{\alpha} ^N) f(t_{\beta} ^ N) = f(t_{\beta} ^ N) f(t_{\alpha}
^ N)$ (since $f$ is injective on K) $\Leftrightarrow$
$t_{f_*(\alpha)} ^P t_{f_*(\beta)}^Q = t_{f_*(\beta)}^Q
t_{f_*(\alpha)} ^P$ where $P = M(\alpha, N), Q = M(\beta, N) \in
\mathbb{Z}^*$ $ \Leftrightarrow i(f_{*}(\alpha),
f_{*}(\beta))=0$.\end{proof}\\

Now, we prove the second main theorem of the section.

\begin{theorem}
\label{main4} Let $K$ be a finite index subgroup of $Mod_R^*$ and
$f$ be an injective homomorphism $f:K \rightarrow Mod_R^*$. If $g
= 2$ and $p=1$ then $f$ has the form $k \rightarrow hkh^{-1}$ for
some $h \in Mod_R^*$ and $f$ has a unique extension to an
automorphism of $Mod_R^*$. If $R$ is a closed surface of genus 2,
then $f$ has the form $k \rightarrow hkh^{-1} i^{m(k)}$ for some
$h \in Mod_R^*$ where $m$ is a homomorphism $K \rightarrow
\mathbb{Z}_2$ and $i$ is the hyperelliptic involution on $R$.
\end{theorem}

\begin{proof} If $g=2$ and $p \leq 1$, by Lemma \ref{intersection0} $f_*$
is a superinjective simplicial map on $\mathcal{C}(R)$. Then by
Theorem \ref{theorem3} $f_*$ is induced by a homeomorphism $h:R
\rightarrow R$, i.e. $f_*(\alpha) = h_\#(\alpha)$ for all $\alpha$
in $\mathcal{A}$, where $h_\#=[h]$. Let $\chi ^ {h\#}: Mod_R^*
\rightarrow Mod_R^*$ be the isomorphism defined by the rule $\chi
^ {h_\#}(k) = h_\#k{h_\#}^{-1}$ for all $k$ in $Mod_R^*$. Then for
all $\alpha$ in $\mathcal{A}$, we have the following:

$\chi ^{h_\# ^{-1}} \circ f ({t_ \alpha} ^N) =  \chi ^{h_\# ^{-1}}
(t_{f_*(\alpha)}^ M) = \chi ^{h_\# ^{-1}} (t_{h_\#(\alpha)} ^M) =
h_\#^{-1} t_{h_\#(\alpha)} ^M h_\# = t_ {h_\#^{-1} (h_\#(\alpha))}
^{M \cdot \epsilon{(h_\#^{-1})}} = t_\alpha ^{M \cdot
\epsilon{(h_\#^{-1})}}$.\\

If $g=2$ and $p=1$, then since $\chi ^{h_\#^{-1}} \circ f$ is
injective, $\chi ^{h_\#^{-1}} \circ f = id_K$ by Lemma
\ref{identity}. So, $\chi ^h_\# |_K = f$. Hence, $f$ is the
restriction of an isomorphism which is conjugation by $h_\#$,
(i.e. $f$ is induced by $h$). Suppose that there exists an
automorphism $\tau :  Mod_R^* \rightarrow Mod_R^*$ such that $\tau
|_{K}=f$. Let $N \in Z^*$ such that $ t_\alpha ^N \in K$  for all
$\alpha$ in $\mathcal{A}$. Since $\chi ^h_\# |_K = f = \tau |_K$
and $t_\alpha ^N \in K$, $\tau(t_\alpha ^N) = \chi ^h_\#(t_\alpha
^N)$ for all $\alpha$ in $\mathcal{A}$. Then by Corollary
\ref{id}, $\tau = \chi ^ {h_\#}$. Hence, the extension of $f$ is unique.\\

If $g=2$ and $p=0$, then since $\chi ^{h_\#^{-1}} \circ f$ is
injective, $\chi ^{h_\#^{-1}} \circ f (k) = ki^{m(k)}$ where $m(k)
\in \{0,1\}$ by Lemma \ref{identity}. So, $f$ has the form $k
\rightarrow h_\# k h_\# ^{-1} i^{m(k)}$. Since $\chi ^{h_\# ^{-1}}
\circ f (k_1 k_2) = k_1 k_2 i^{m(k_1 k_2)}$ and $\chi ^{h_\#
^{-1}} \circ f (k_1) \chi ^{h^{-1}_\# } \circ f ( k_2)$ $ = k_1
i^{m(k_1)} k_2 i^{m(k_2)} = k_1 k_2 i^{m(k_1)} i^{m(k_2)}$, we
have that $m(k_1 k_2) = m(k_1) + m(k_2)$ for all $k_1, k_2 \in K$.
So, $m: K \rightarrow \mathbb{Z}_2$ is a homomorphism.\end{proof}\\

\noindent {\bf Remark:} Note that $k \rightarrow hkh^{-1}
i^{m(k)}$ defines a homomorphism from $K \rightarrow Mod_R^*$ for
every $h \in Mod_R^*$ and for every homomorphism $m : K
\rightarrow \mathbb{Z}_2$. It is easy to see that $k \rightarrow
hkh^{-1} i^{m(k)}$ is injective if and only if either $i \notin K$
or $i \in Ker(m)$. Inner automorphisms of $K$ act on the set of
injective homomorphisms from $K \rightarrow Mod_R^*$. By using
Theorem \ref{main4} we can see that the orbit space $InjHom(K,
Mod_R^*)/Inn(K)$ of this action is finite. Then we have the
following corollary.

\begin{coroll} Suppose that $g=2$ and $p \leq 1$. Let $K$ be a finite
index subgroup of $Mod_R^*$. Then $Out(K)$ is finite.
\end{coroll}

In the other cases, when $R$ has genus at least two and $K$ is a
finite index subgroup, we have that $Out(K)$ is finite as a
corollary to the main results in \cite{Ir1}, \cite{Ir2}. See
\cite{Mc1} for an explicit description of automorphisms of
$Mod_R^*$ for a closed surface of genus two.

\section{Extending Superinjective Simplicial Maps of $\mathcal{N}(R)$ to
Superinjective Simplicial Maps of $\mathcal{C}(R)$}

In this section we prove Theorem \ref{theorem1} and Theorem
\ref{theorem2}.

\begin{lemma} \label{extension} Suppose that $g \geq 2$ and $p \leq 1$.
Let $\lambda : \mathcal{N}(R) \rightarrow
\mathcal{N}(R)$ be a superinjective simplicial map. Then $\lambda$
extends to a superinjective simplicial map $\lambda_* :
\mathcal{C}(R) \rightarrow \mathcal{C}(R)$.\end{lemma}

\begin{proof} If $x$ is a nonseparating simple closed curve, we define
$\lambda_*([x])=\lambda([x])$. Let $c$ be a separating simple
closed curve on $R$.\\

Case 1: Assume that $R$ is closed. Since $g \geq 2$, $c$ separates
$R$ into two subsurfaces $R_1, R_2$, and both of $R_1, R_2$ have
genus at least one. We take a chain on $R_1$, $\{\alpha_1, ...,
\alpha_m\}$ with $i( \alpha_{i}, \alpha_{i+1})=1$, $i( \alpha_{i},
\alpha_{j})=0$ for $|i-j| > 1$, $[a_i] = \alpha_i \in
\mathcal{N}(R)$, as shown in Figure 8, i (for $g=4$ case when
$R_1$ has genus 3), such that $R_1 \cup \{c\}$ is a regular
neighborhood of $a_1 \cup ... \cup a_n$. Since $\lambda$ preserves
disjointness, nondisjointness and intersection one property, we
can see that the chain $\{\alpha_1, ..., \alpha_n \}$ is mapped by
$\lambda$ into a similar chain, $\{\lambda(\alpha_1), ...,
\lambda(\alpha_m) \}$ with $i( \lambda(\alpha_{i}),
\lambda(\alpha_{i+1}))=1$, $i(\lambda(\alpha_{i}),
\lambda(\alpha_{j}))=0$ for $|i-j| > 1$. Let $a_i' \in
\lambda(\alpha_i)$ such that any two elements in $\{a_1', ...,
a_m'\}$ have minimal intersection with each other. Let $M$ be a
regular neighborhood of $a_1' \cup ... \cup a_n'$. Then it is easy
to see that $M$ is homeomorphic to $R_1 \cup c$. Let $a'$ be the
boundary of $M$. Suppose that we have another chain on $R_1$,
$\{\beta_1, ..., \beta_m \}$ with $i( \beta_{i}, \beta_{i+1})=1$,
$i( \beta_{i}, \beta_{j})=0$ for $|i-j| > 1$, $[b_i] = \beta_i \in
\mathcal{N}(R)$ such that $R_1 \cup \{c\}$ is a regular
neighborhood of $b_1 \cup ... \cup b_m$. Again we see that the
chain $\{\beta_1, ..., \beta_m \}$ is mapped by $\lambda$ into a
similar chain, $\{\lambda(\beta_1), ..., \lambda(\beta_m) \}$ with
$i( \lambda(\beta_{i}), \lambda(\beta_{i+1}))=1$,
$i(\lambda(\beta_{i}), \lambda(\beta_{j}))=0$ for $|i-j| > 1$. Let
$b_i' \in \lambda(\beta_i)$ such that any two elements in $\{b_1',
..., b_m'\}$ have minimal intersection with each other. Let $T$ be
a regular neighborhood of $b_1' \cup ... \cup b_n'$. Then $T$ is
homeomorphic to $R_1 \cup c$. Let $b'$ be the boundary of $T$.\\

Claim: $[a'] = [b']$.\\

Proof: We take a similar chain on $R_2$, $\{\gamma_1, ...,
\gamma_n \}$ with $i(\gamma_{i}, \gamma_{i+1})=1$, $i( \gamma_{i},
\gamma_{j})=0$ for $|i-j| > 1$, $[c_i] = \gamma_i \in
\mathcal{N}(R)$, such that $R_2 \cup \{c\}$ is a regular
neighborhood of $c_1 \cup ... \cup c_n$. This chain is mapped into
a similar chain, $\{\lambda(\gamma_1), ..., \lambda(\gamma_n) \}$.
Let $c_i' \in \lambda(\gamma_i)$ such that any two elements in
$\{a_1', ..., a_n', c_1', ..., c_m'\}$ and $\{b_1', ..., b_n',
c_1', ..., c_m'\}$ have minimal intersection with each other. Then
$a_i'$ is disjoint from $c_j'$ for any $i, j$, and $b_i'$ is
disjoint from $c_j'$ for any $i, j$. Then we can choose a regular
neighborhood $N$ of $c_1' \cup ... \cup c_n'$ in $R$ such that $N$
is disjoint from $M \cup T$. Let $c'$ be the boundary component of
$N$. Then it is easy to see that $N$ is homeomorphic to $R_2 \cup
c$ and the boundary components of $M$ and $N$ are isotopic in $R$.
Similarly, the boundary components of $T$ and $N$ are isotopic in
$R$. Hence, $[a'] = [c'] = [b']$. We define $\lambda_*([c]) = [a']$.\\

Claim: $\lambda_*: \mathcal{C}(R) \rightarrow \mathcal{C}(R)$ is
a simplicial map.\\

Proof: Let $\alpha, \beta$ be two vertices in $\mathcal{C(R)}$
such that $i(\alpha, \beta)=0$. Let $x$ and $y$ be disjoint
representatives of $\alpha$ and $\beta$ respectively. If $x, y$
are two nonseparating simple closed curves then we have
$i(\lambda_*(\alpha), \lambda_*(\beta)) = i(\lambda(\alpha),
\lambda(\beta))=0$. If $x$ is a nonseparating simple closed curve
and $y$ is a separating simple closed curve then $x$ lives in a
subsurface $R_1$ which comes from separation by $y$. Then we could
choose a chain as described above which contains $x$, and then see
that by the construction of the image of $[y]$, we get that
$i(\lambda_*([x]), \lambda_*([y]))=0$. If both $x$ and $y$ are
separating, then it is easy to see that $R \setminus ((R \setminus
{x}) \cap (R \setminus {y}))$ has two connected components $T_1,
T_2$ such that $T_1$ is disjoint from $T_2$ and $y$ is an
essential boundary component of $T_1$ and $x$ is an essential
boundary component of $T_2$. Then we see that the chains which
come from disjoint subsurfaces $T_1$ and $T_2$ which will be used
to define the images of $[x]$ and $[y]$ are disjoint. Since,
$\lambda$ preserves disjointness we see that the ``image" chains
will be disjoint, and hence, $i(\lambda_*([x]),
\lambda_*([y]))=0$.\\

\begin{figure}
\begin{center}
\epsfxsize=4.2in \epsfbox{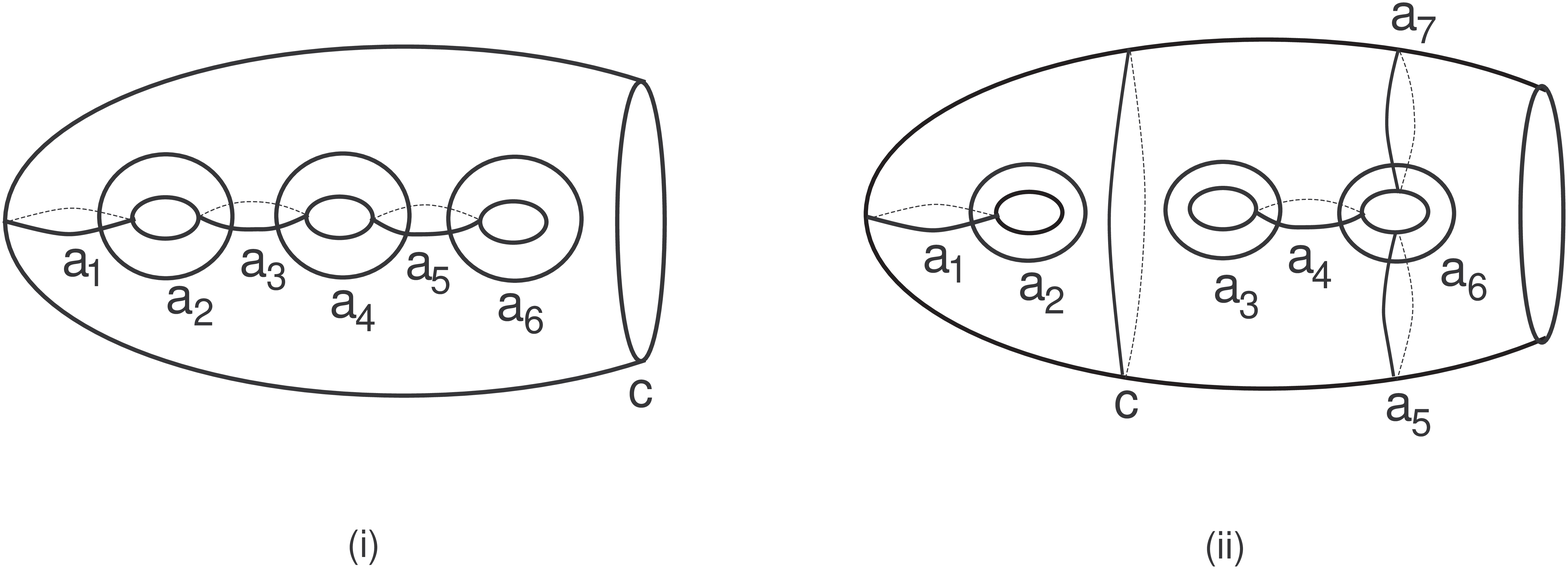} \caption{Chains}
\end{center}
\end{figure}

Claim: $\lambda_*: \mathcal{C}(R) \rightarrow \mathcal{C}(R)$ is
a superinjective simplicial map.\\

Proof: Let $\alpha, \beta$ be two vertices in $\mathcal{C(R)}$
such that $i(\alpha, \beta) \neq 0$. Let $x$ and $y$ be
representatives of $\alpha$ and $\beta$ respectively.\\

(i) Assume that $x$ and $y$ are two nonseparating simple closed
curves. Then we have $i(\lambda_*(\alpha), \lambda_*(\beta)) =
i(\lambda(\alpha), \lambda(\beta)) \neq 0$ (since $\lambda$ is
superinjective).\\

(ii) Assume that $y$ is a nonseparating simple closed curve and
$x$ is a separating simple closed curve, and $x$ separates $R$
into two subsurfaces $R_1, R_2$. W.L.O.G. assume that $y$ is in
$R_1$. Let $\{\alpha_1, ..., \alpha_m\}$ be a chain on $R_1$ with
$i( \alpha_{i}, \alpha_{i+1})=1$, $i( \alpha_{i}, \alpha_{j})=0$
for $|i-j| > 1$, $[a_i] = \alpha_i \in \mathcal{N}(R)$ such that
$R_1 \cup \{x\}$ is a regular neighborhood of $a_1 \cup ... \cup
a_m$. Then, since $i(\alpha, \beta) \neq 0$, $i(\beta, \alpha_i)
\neq 0$ for some $i$. Then $i(\lambda_*(\beta),
\lambda_*(\alpha_i)) = i(\lambda(\beta), \lambda(\alpha_i)) \neq
0$ (since $\lambda$ is superinjective). Then it is easy to see
that $i(\lambda_*(\alpha),
\lambda_*(\beta)) \neq 0$.\\

(iii) Assume that both $x$ and $y$ are separating and $x$
separates $R$ into two subsurfaces $R_1, R_2$. Let $\{\alpha_1,
..., \alpha_m\}$ be a chain on $R_1$ with $i( \alpha_{i},
\alpha_{i+1})=1$, $i( \alpha_{i}, \alpha_{j})=0$ for $|i-j| > 1$,
$[a_i] = \alpha_i \in \mathcal{N}(R)$ such that $R_1 \cup \{x\}$
is a regular neighborhood of $a_1 \cup ... \cup a_m$. Then, since
$i(\alpha, \beta) \neq 0$, $i(\beta, \alpha_i) \neq 0$ for some
$i$. Then $i(\lambda_*(\beta), \lambda_*(\alpha_i)) \neq 0$ by
(ii). Then it is easy to see that $i(\lambda_*(\alpha),
\lambda_*(\beta)) \neq 0$. Hence, we have a superinjective
simplicial extension $\lambda_*: \mathcal{C}(R) \rightarrow
\mathcal{C}(R)$ of $\lambda: \mathcal{N}(R) \rightarrow
\mathcal{N}(R)$.\\

Case 2: Assume that $R$ has one boundary component. Then $c$
separates $R$ into two subsurfaces $R_1, R_2$. W.L.O.G. assume
that $R_1$ is a genus $k$ subsurface having $c$ as its boundary.
We consider chains on $R_1$ as in the first case, and chains on
$R_2$ such that regular neighborhoods of the curves coming from
the chains have $c$ as their essential boundary component and the
boundary of $R$ as their inessential boundary component (see
Figure 7, ii). Then it is easy to see that the proof of the lemma
is similar to case 1.\end{proof}\\

\noindent {\bf Remark:} If $g \geq 2$ and $p \geq 2$ and $C$ is
the set of separating circles on $R$ which separate $R$ into two
pieces such that each piece has genus at least one, then by using
chains on these two pieces and following the techniques in the
previous lemma we can extend $\lambda$ to $\lambda_*$ on $C$ and
get a superinjective extension.\\

Let $M$ be a sphere with $k$ holes and $k \geq 5$. A circle $a$ on
$M$ is called an {\it n-circle} if $a$ bounds a disk with $n$
holes on $M$ where $n \geq 2$. A {\it pentagon} in
$\mathcal{C}(M)$ is an ordered 5-tuple $(\alpha_1, \alpha_2,
\alpha_3, \alpha_4, \alpha_5)$, defined up to cyclic permutations,
of vertices of $\mathcal{C}(M)$ such that $i(\alpha_j,
\alpha_{j+1}) = 0$ for $j=1,2,...,5$ and $i(\alpha_j, \alpha_k)
\neq 0$ otherwise, where $\alpha_6 = \alpha_1$. A vertex in
$\mathcal{C}(M)$ is called an {\it n-vertex} if it has a
representative which is an n-circle on $M$.\\

Let $x, y$ be disjoint simple closed curves on $R$ such that
$([x], [y])$ is a peripheral pair, i.e. $x, y$ and a boundary
component of $R$ bound a pair of pants on $R$. $x \cup y$ separate
$R$ into two subsurfaces. Let $R_{x, y}$ be the positive genus
subsurface of $R$ which comes from this separation. We can
identify $\mathcal{N}(R_{x, y})$ with a subcomplex, $L_{x,y}$ of
$\mathcal{N}(R)$. By $\lambda_{x, y}$ we will denote the
restriction of $\lambda$ on $\mathcal{N}(R_{x, y})$. If $k: R
\rightarrow R$ is a homeomorphism then we will use $k_\#$ for the
map induced by $k$ on $\mathcal{N}(R)$ (i.e. $k_\#: \mathcal{N}(R)
\rightarrow \mathcal{N}(R)$ where $k_\# =[k]$).\\

An arc $i$ on $R$ is called \textit{properly embedded} if
$\partial i \subseteq \partial R$ and $i$ is transversal to
$\partial R$. $i$ is called \textit{nontrivial} (or
\textit{essential}) if $i$ cannot be deformed into $\partial R$ in
such a way that the endpoints of $i$ stay in $\partial R$ during
the deformation. If $a$ and $b$ are two disjoint arcs connecting a
boundary component of $R$ to itself, $a$ and $b$ are called {\it
linked} if their end points alternate on the boundary component.
Otherwise, they are called {\it unlinked}. The \textit{complex of
arcs}, $\mathcal{B}(R)$, on $R$ is an abstract simplicial complex.
Its vertices are the isotopy classes of nontrivial properly
embedded arcs $i$ in $R$. A set of vertices forms a simplex if
these vertices can be represented by pairwise disjoint arcs. Let
$i$ be an essential properly embedded arc on $R$. Let $A$ be a
boundary component of $R$ which has one end point of $i$ and $B$
be the boundary component of $R$ which has the other end point of
$i$. Let $N$ be a regular neighborhood of $i \cup A \cup B$ in
$R$. By Euler characteristic arguments, $N$ is a pair of pants.
The boundary components of $N$ are called {\it encoding circles of
$i$ on $R$}. An essential properly embedded arc $i$ on $R$ is
called {\it type 1} if it joins one boundary component
$\partial_k$ of $R$ to itself. $i$ is called {\it nonseparating}
if its complement in $R$ is connected.\\

The mapping class group, $Mod_R$, of $R$ is the group of isotopy
classes of orientation preserving homeomorphisms of $R$. The pure
mapping class group, $PMod_R$, is the subgroup of $Mod_R$
consisting of isotopy classes of homeomorphisms which preserve
each boundary component of $R$.

\begin{lemma}
\label{imp} Suppose $g \geq 3$ and $p=2$. Let $\lambda:
\mathcal{N}(R) \rightarrow \mathcal{N}(R)$ be a superinjective
simplicial map. Assume that for any peripheral pair $([x], [y])$
on $R$ with $x, y$ disjoint, $\lambda_{x, y}$ agrees with a map,
$(g_{x, y})_\# : \mathcal{N}(R_{x, y}) \rightarrow
\mathcal{N}(R_{x', y'})$, which is induced by a homeomorphism
$g_{x, y} : R_{x,y} \rightarrow R_{x', y'}$ where $x', y'$ are
disjoint, $\lambda([x]) = [x']$, $\lambda([y]) = [y']$, $g_{x,
y}(x)= x'$, $g_{x, y}(y)= y'$. Then $\lambda$ agrees with a map
$h_\#: \mathcal{N}(R) \rightarrow \mathcal{N}(R)$ which is induced
by a homeomorphism $h: R \rightarrow R$.\end{lemma}

\begin{figure}
\begin{center}
\epsfxsize=4.7in \epsfbox{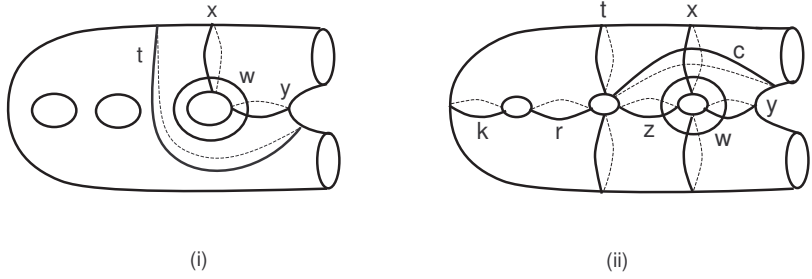} \caption{(i) A dual curve,
(ii) Pants decomposition}
\end{center}
\end{figure}

\begin{proof} Let $x, y$ be disjoint simple closed curves such
that $([x], [y])$ is a peripheral pair, and let $(g_{x, y})_\# :
\mathcal{N}(R_{x, y}) \rightarrow \mathcal{N}(R_{x', y'})$ be a
simplicial map which is induced by a homeomorphism $g_{x, y} :
R_{x,y} \rightarrow R_{x', y'}$ where $x', y'$ are disjoint,
$\lambda([x]) = [x']$, $\lambda([y]) = [y']$, $g_{x, y}(x)= x'$,
$g_{x, y}(y)= y'$ such that $\lambda_{x, y}$ agrees with $(g_{x,
y})_\# $ on $\mathcal{N}(R_{x, y})$. Let $g$ be a homeomorphism of
$R$ which cuts to a homeomorphism $R_{x,y} \rightarrow R_{x', y'}$
which is isotopic to $g_{x, y}$. Then each homeomorphism of $R$
which cuts to a homeomorphism $R_{x, y} \rightarrow R_{x', y'}$
which is isotopic to $g_{x, y}$, is isotopic to an element in the
set $\{gt_x^mt_y^n, m, n \in \mathbb{Z}\}$. It is easy to see that
$\lambda_{x, y}$ agrees with the restriction of $(gt_x^mt_y^n)_\#$
on $\mathcal{N}(R_{x, y})$ for all $m, n \in \mathbb{Z}$. Let $w$
be a simple closed curve which is dual to both of $x$ and $y$ (see
Figure 9 (i)). Since $\lambda$ preserves geometric intersection
one property by Lemma \ref{intone}, $\lambda([w])$ has a
representative which is dual to both of $x'$ and $y'$. Let $P$ be
a regular neighborhood of $x \cup y \cup w$. $P$ is a genus one
surface with two boundary components. Let $t$ be the essential
boundary component of $P$. Let $Q$ be the genus one subsurface
with two boundary components of $R$ which has $g(t)$ as its
boundary. Then by using the properties of $\lambda$ it is easy to
see that $\lambda([w])$ has a representative which lies in the
interior of $Q$, and dual to both of $x'$ and $y'$ and there
exists $m_o, n_o \in \mathbb{Z}$ such that $gt_x^{m_o} t_y^{n_o}$
agrees with $\lambda$ on $[w]$. Let $D_{x,y}$ be the set of
isotopy classes of simple
closed curves which are dual to each of $x$ and $y$ on $R$.\\

Claim 1: $(gt_x^{m_o} t_y^{n_o})_\#$ agrees with $\lambda$ on
$\{[x]\} \cup \{[y]\} \cup L_{x,y} \cup D_{x,y}$.\\

Proof: It is clear that $(gt_x^{m_o} t_y^{n_o})_\# ([x])= \lambda
([x]) = [x']$ and $(gt_x^{m_o} t_y^{n_o})_\# ([y])= \lambda ([y])
= [y']$. Since $\lambda_{x, y}$ agrees with the restriction of
$(gt_x^{m_o} t_y^{n_o})_\#$ on $\mathcal{N}(R_{x, y})$,
$(gt_x^{m_o} t_y^{n_o})_\#$ agrees with $\lambda$ on $L_{x, y}$.
We have seen that $(gt_x^{m_o} t_y^{n_o})_\#$ agrees with
$\lambda$ on $[w]$. Let $w_1$ be a simple closed curve which is
disjoint from $w$ and dual to both of $x$ and $y$. As we described
before, there exists $\tilde{m}, \tilde{n} \in \mathbb{Z}$ such
that $\lambda$ agrees with $(gt_x^{\tilde{m}} t_y^{\tilde{n}})_\#$
on $[w_1]$. Since $w$ and $w_1$ are disjoint nonseparating curves,
$i(\lambda([w]), \lambda([w_1])) = 0$. Then by using the
properties that $\lambda$ preserves disjointness and
nondisjointness, we can see that $m = \tilde{m}$ and $n =
\tilde{n}$. This shows that $(gt_x^{m_o} t_y^{n_o})_\#$ also
agrees with $\lambda$ on $[w_1]$. Given any simple closed curve
$v$ which is dual to both of $x$ and $y$, we can find a sequence
of dual curves to both of $x$ and $y$, connecting $w$ to $v$, such
that each consecutive pair is disjoint, i.e. the isotopy classes
of these curves define a path between $w$ and $v$ in
$\mathcal{N}(R)$. Then using the argument given above and the
sequence, we conclude that $(gt_x^{m_o} t_y^{n_o})_\#$ agrees with
$\lambda$ on $D_{x, y}$. Hence, $(gt_x^{m_o} t_y^{n_o})_\#$ agrees
with $\lambda$ on $\{[x]\} \cup \{[y]\} \cup L_{x, y} \cup D_{x,
y}$. This proves claim 1. Let $h_{x,y} = gt_x^{m_o} t_y^{n_o}$. We
have that $(h_{x,y})_\#$ agrees with $\lambda$ on
$\{[x]\} \cup \{[y]\} \cup L_{x, y} \cup D_{x, y}$.\\

We complete $x, y$ to a pair of pants decomposition $P$ consisting
of nonseparating circles as shown in Figure 9 (ii) for $g=3, p=2$,
(similar configurations can be chosen for $g \geq 3$).
Let $t, z, w, c, k, r$ be as shown in the figure.\\

\begin{figure}
\begin{center}
\epsfxsize=2.55in \epsfbox{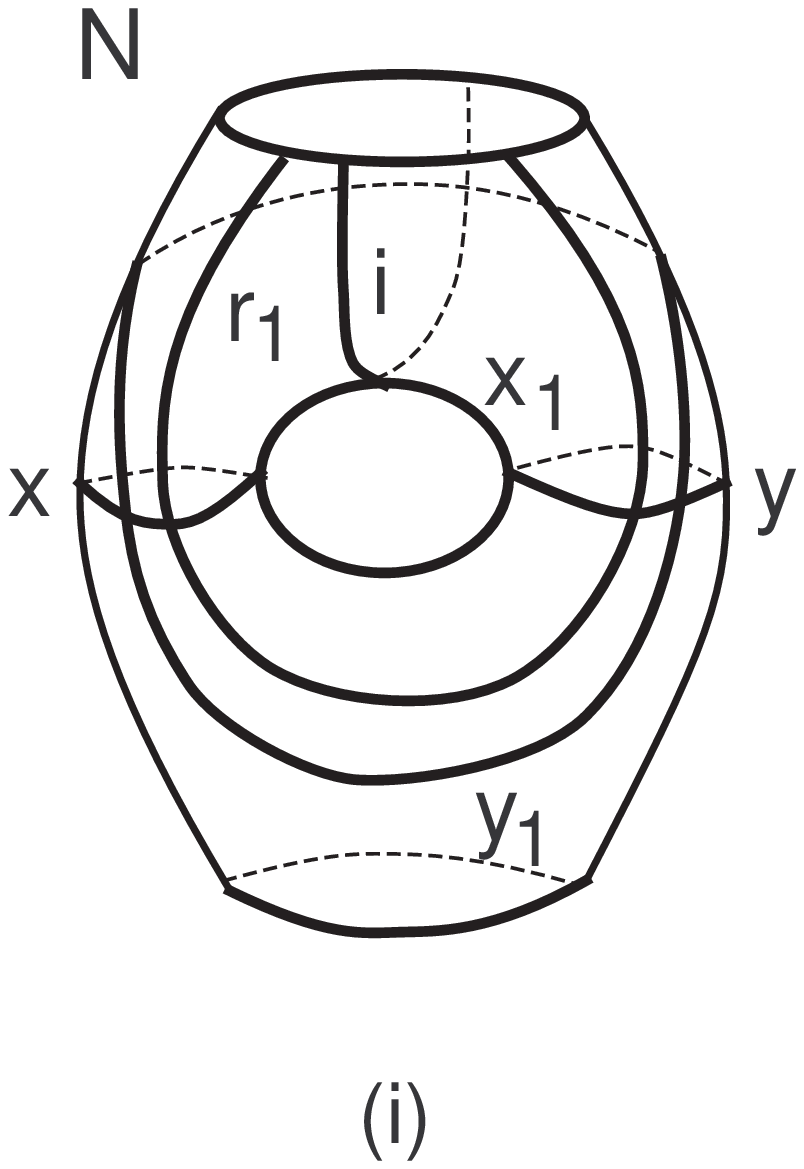} \epsfxsize=2in
\epsfbox{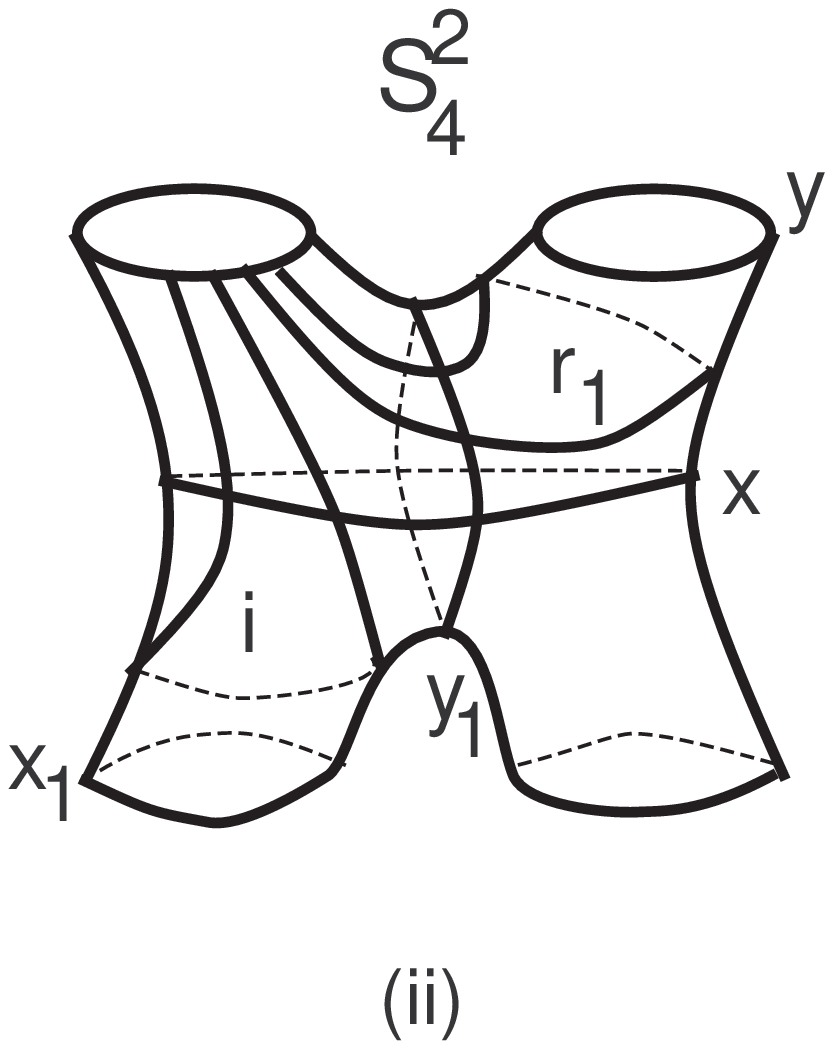} \caption{Arcs and their encoding circles}
\end{center}
\end{figure}

Claim 2: $(h_{x,y})_\#$ agrees with $\lambda$ on $c$.\\

Proof: Let $t', z', x', y', k', r'$ be pairwise disjoint
representatives of $\lambda([t]), \lambda([z]), \lambda([x]),$
$\lambda([y]), \lambda([k]), \lambda([r])$ respectively and let
$w'$ be a representative of $\lambda([w])$ which has minimal
intersection with each of $t', z', x', y', k', r'$. Since $k, r,
z, y, \partial_1$ bound a sphere with five holes, $T$, on $R$,
containing $x$ and $c$, then $k', r', z', y'$ and a boundary
component of $R$ bound a sphere with five holes, $T'$, on $R$ by
Lemma \ref{top}. Since $x$ and $c$ have geometric intersection 2
and algebraic intersection 0 in $T$, it is easy to see that there
exist vertices $\gamma_1, \gamma_2, \gamma_3$ of $\mathcal{C}(T)$
such that $(\gamma_1, \gamma_2, [x], \gamma_3, [c])$ is a pentagon
in $\mathcal{C}(T)$, $\gamma_1$ and $\gamma_3$ are 2-vertices,
$\gamma_2$ is a 3-vertex, and $\{[x], \gamma_3\}$, $\{[x],
\gamma_2\}$, $\{[c], \gamma_3\}$ and $\{\gamma_1, \gamma_2 \}$ are
codimension-zero simplices of $\mathcal{C}(T)$, and each of
$\gamma_i$ has a representative which is nonseparating on $R$.
Since $\lambda$ is superinjective, we can see that
$(\lambda(\gamma_1), \lambda(\gamma_2)$, $\lambda([x]),
\lambda(\gamma_3)$, $\lambda([c]))$ is a pentagon in
$\mathcal{C}(T')$. By Lemma \ref{top}, $\lambda(\gamma_1)$ and
$\lambda(\gamma_3)$ are 2-vertices, and $\lambda(\gamma_2)$ is a
3-vertex in $\mathcal{C}(T')$. Since $\lambda$ is an injective
simplicial map $\{\lambda([x]), \lambda(\gamma_3)\}$,
$\{\lambda([x]), \lambda(\gamma_2)\}$, $\{\lambda([c])$,
$\lambda(\gamma_3)\}$ and $\{\lambda(\gamma_1),
\lambda(\gamma_2)\}$ are codimension-zero simplices of
$\mathcal{C}(T')$. It is easy to see that $\lambda([c])$ has a
representative which is disjoint from $t', z', y'$. Then
$\{\lambda([x]), \lambda([c])\}$ have representatives with
geometric intersection 2 and algebraic intersection 0 in the
sphere with four holes bounded by $t', z', y'$ and the boundary
component of $R$, \cite{K}. We also have that $(h_{x,y})_\# ([c])$
has a representative $c''$ which is disjoint from $t', z', y', w'$
and which intersects $x'$ geometrically twice and algebraically 0
times. Since $(h_{x,y})_\#$ agrees with $\lambda$ on $\{[x]\} \cup
\{[y]\} \cup L_{x, y} \cup D_{x, y}$ it is easy to see that
$[c']=[c'']$, i.e.
$(h_{x,y})_\#$ agrees with $\lambda$ on $c$.\\

Claim 3: $(h_{x,y})_\#$ agrees with $\lambda$ on the class of
every nonseparating circle on $R$.\\

Proof: Let $z$ be a nonseparating simple closed curve on $R$. Let
$t$ be another simple closed curve disjoint from $z$ such that $z,
t$ and $\partial_1$ bound a pair of pants in $R$. Let $i$ and $j$
be nonseparating type 1 arcs connecting $\partial_1$ to itself
such that $i$ has $x, y$ as its encoding circles and $j$ has $z,
t$ as its encoding circles. W.L.O.G. we can assume that $i$ and
$j$ have minimal intersection. By Lemma 3.8 in \cite{Ir2}, there
is a sequence $i = r_0 \rightarrow r_1 \rightarrow ... \rightarrow
r_{n+1}=j$ of essential properly embedded nonseparating type 1
arcs joining $\partial_1$ to itself so that each consecutive pair
is disjoint, i.e. the isotopy classes of these arcs define a path
in $\mathcal{B}(R)$ between $i$ and $j$. Let $x_i, y_i$ be the
encoding circles for $r_i$ for $i=1,..., n$. For the pair of arcs
$i$ and $r_1$, we will consider the following cases:\\

Case i: Assume that $i$ and $r_1$ are linked, i.e their end points
alternate on the boundary component $\partial_1$. Then a regular
neighborhood of $i \cup r_1 \cup \partial_1$ is a genus one
surface with two boundary components $N$, and the arcs $i, r_1$
and their encoding circles $x, y, x_1, y_1$ on $N$ are as shown in
Figure 10 (i). In this case, we complete $\{x, y, x_1, y_1\}$ to a
curve configuration $G$ consisting of nonseparating circles which
is shown in Figure 11 (i), for $g=3, p=3$, see \cite{IMc}, such
that the isotopy classes of Dehn twists about the elements of this
set generate $PMod_R$ and all the curves in this set are (i)
either disjoint from $x, y$ or simultaneous dual to $x, y$, and
(ii) either disjoint from $x_1, y_1$ or simultaneous dual to $x_1,
y_1$. Then, since all the curves in $G$ are either disjoint from
$x, y$ or simultaneous dual to $x, y$, by claim 1 we have that
$(h_{x,y})_\# ([x]) = \lambda ([x])$ for every $x \in G$.
Similarly, since all the curves in $G$ are either disjoint from
$x_1, y_1$ or simultaneous dual to $x_1, y_1$, by claim 1 we have
$(h_{x_1,y_1})_\# ([x]) = \lambda ([x])$ for every $x \in G$.
Hence $(h_{x,y})_\# ([x])= \lambda ([x]) = (h_{x_1,y_1})_\# ([x])$
for every $x \in G$. Then $(h_{x,y}^{-1} h_{x_1,y_1})_\# \in
C_{Mod_R}(PMod_R)$. By Theorem 5.3 in \cite{IMc},
$C_{Mod_R}(PMod_R) = \{1\}$. Hence $(h_{x,y})_\# = (h_{x_1,y_1})_\#$.\\

\begin{figure}
\begin{center}
\epsfxsize=2.35in \epsfbox{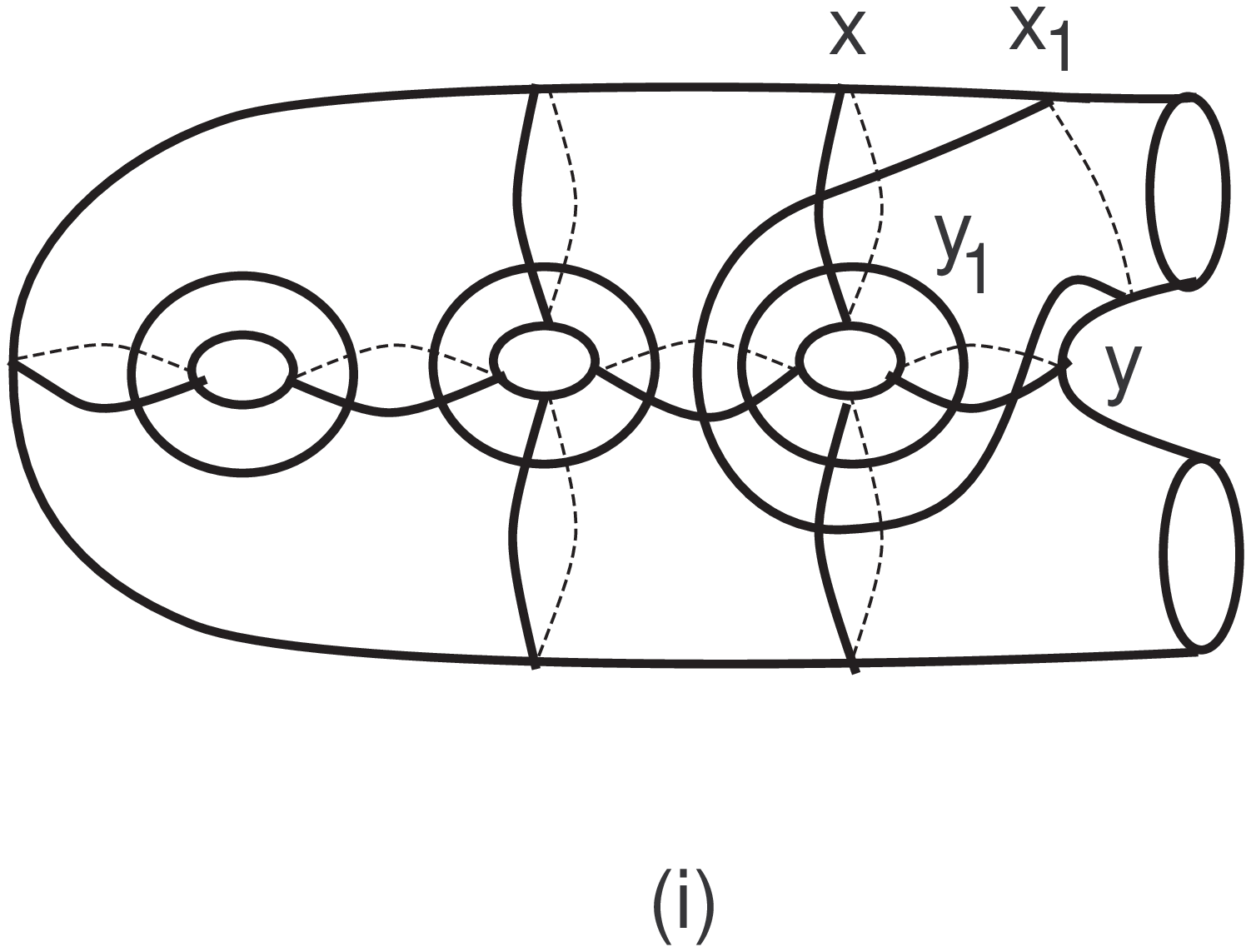} \epsfxsize=2.35in
\epsfbox{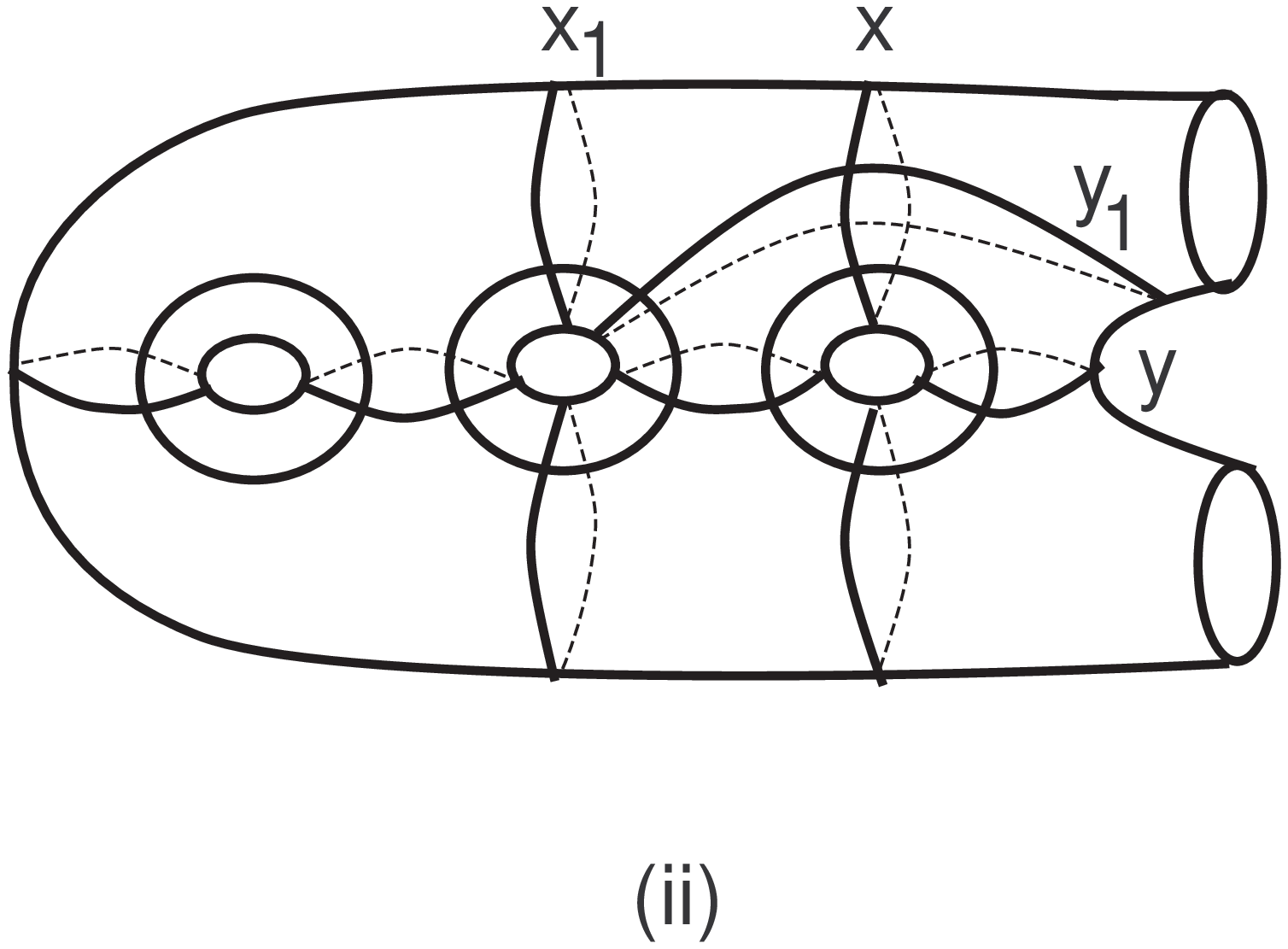} \caption{Encoding circles in a curve
configuration}
\end{center}
\end{figure}

Case ii: Assume that $i$ and $r_1$ are unlinked, i.e their end
points don't alternate on the boundary component $\partial_1$.
Then a regular neighborhood of $i \cup r_1 \cup \partial_1$ is a
sphere with four boundary components $S_4 ^2$, and the arcs $i,
r_1$ and their encoding circles $x, y, x_1, y_1$ on $S_4 ^2$ are
as shown in Figure 10 (ii). Let $w$ be the boundary component of
$S_4 ^2$ which is different from $x_1, y, \partial_1$. If $w$ is a
nonseparating curve, then we complete $\{x, y, x_1, y_1\}$ to a
curve configuration $G$ consisting of nonseparating circles as
shown in Figure 11 (ii) such that the isotopy classes of Dehn
twists about the elements of $G$ generate $PMod_R$. By claim 1 and
claim 2, $(h_{x,y})_\# ([x]) = \lambda ([x])$ for every $x \in G$
and $(h_{x_1,y_1})_\# ([x]) = \lambda ([x])$ for every $x \in G$.
Then $(h_{x,y}^{-1} h_{x_1,y_1})_\# ([x]) = [x]$ for every $x \in
G$. Then $(h_{x,y}^{-1} h_{x_1,y_1})_\# \in C_{Mod_R}(PMod_R)$. By
Theorem 5.3 in \cite{IMc}, $C_{Mod_R}(PMod_R) = \{1\}$. Hence
$(h_{x,y})_\#= (h_{x_1,y_1})_\#$.\\

Suppose that $w$ is a separating curve. By the remark after Lemma
\ref{extension}, we can extend $\lambda$ to a superinjective map
$\lambda_*$ on a subcomplex of $\mathcal{C}(R)$ containing
separating circles on $R$ which separate $R$ into two pieces such
that each piece has genus at least one. Notice that $w$ is such a
circle. For the rest of the proof, we will use $\lambda$ for this
extension. We will do the case when $x_1$ and $y$ are not
isotopic. The other case can be done similarly. Let $M$ be the
subsurface which has $w$ on its boundary and which does not
contain $N$. Let $T$ be the closure of $R \setminus \{M \cup N\}$.
The circles $x_1$ and $y$ are boundary components of $T$. Since
$w$ is an essential separating circle and $p=2$, $M$ has genus at
least one. In Figure 12, we show $M, N, T$ for a special case. By
claim 1, $(h_{x,y})_\# ([x]) = \lambda ([x])$ for every $x \in
\mathcal{N}(T)$. Similarly, $(h_{x_1,y_1})_\# ([x]) = \lambda
([x])$ for every $x \in \mathcal{N}(T)$. Then $(h_{x,y}^{-1}
h_{x_1,y_1})_\# ([x]) = [x]$ for every $x \in \mathcal{N}(T)$.
Then the restriction of $(h_{x,y}^{-1} h_{x_1,y_1})_\#$ on
$\mathcal{N}(T)$ is in $C(PMod_T)$. By Theorem 5.3 in \cite{IMc},
$C(PMod_T) = \{1\}$. Hence $(h_{x,y})_\# = (h_{x_1,y_1})_\#$ on
$\mathcal{N}(T)$. It is easy to see that $(h_{x,y})_\# =
(h_{x_1,y_1})_\#$ on the set $\{[x_1], [w], [y]\}$. Following the
proof of claim 2 (considering that we have the extended
superinjective simplicial map on ``good" separating circles), we
see that $(h_{x,y})_\# = (h_{x_1,y_1})_\#$ on $\{[x], [y_1]\}$.
Then, since Dehn twists about $x$ and $y_1$ generate $PMod_N$, the
restriction of $(h_{x,y}^{-1} h_{x_1,y_1})_\#$ on $\mathcal{C}(N)$
is in $C(PMod_N)$. By Theorem 5.3 in \cite{IMc}, $C(PMod_N) =
\{1\}$. Hence $(h_{x,y})_\# = (h_{x_1,y_1})_\#$ on
$\mathcal{C}(N)$. By claim 1, $(h_{x,y}^{-1} h_{x_1,y_1})_\# ([x])
= [x]$ for every $x \in \mathcal{N}(M)$. Then the restriction of
$(h_{x,y}^{-1} h_{x_1,y_1})_\#$ on $\mathcal{N}(M)$ is in
$C(PMod_M)$. By considering the action on oriented circles and
using Theorem 5.3 in \cite{IMc}, we see that $(h_{x,y})_\# =
(h_{x_1,y_1})_\#$ on $\mathcal{N}(M)$. So we have $(h_{x,y})_\# =
(h_{x_1,y_1})_\#$ on $\mathcal{N}(M) \cup \mathcal{C}(N) \cup
\mathcal{N}(T)$. In Figure 12 we see the curve $c$ which is dual
to each of $x, y, x_1, y_1$. By claim 1, $(h_{x,y})_\# ([c]) =
\lambda ([c]) = (h_{x_1,y_1})_\# ([c])$. Then, since
$(h_{x,y})_\#$ and $(h_{x_1,y_1})_\#$ agree on $[c]$, we have
$(h_{x,y})_\# = (h_{x_1,y_1})_\#$ on $\mathcal{N}(R)$. Then, since
$g \geq 3$, $(h_{x,y})_\# = (h_{x_1,y_1})_\#$. If $w$ is a
boundary component of $R$, then showing that $(h_{x,y})_\# =
(h_{x_1,y_1})_\#$
is similar.\\

\begin{figure}
\begin{center}
\epsfxsize=3.1in \epsfbox{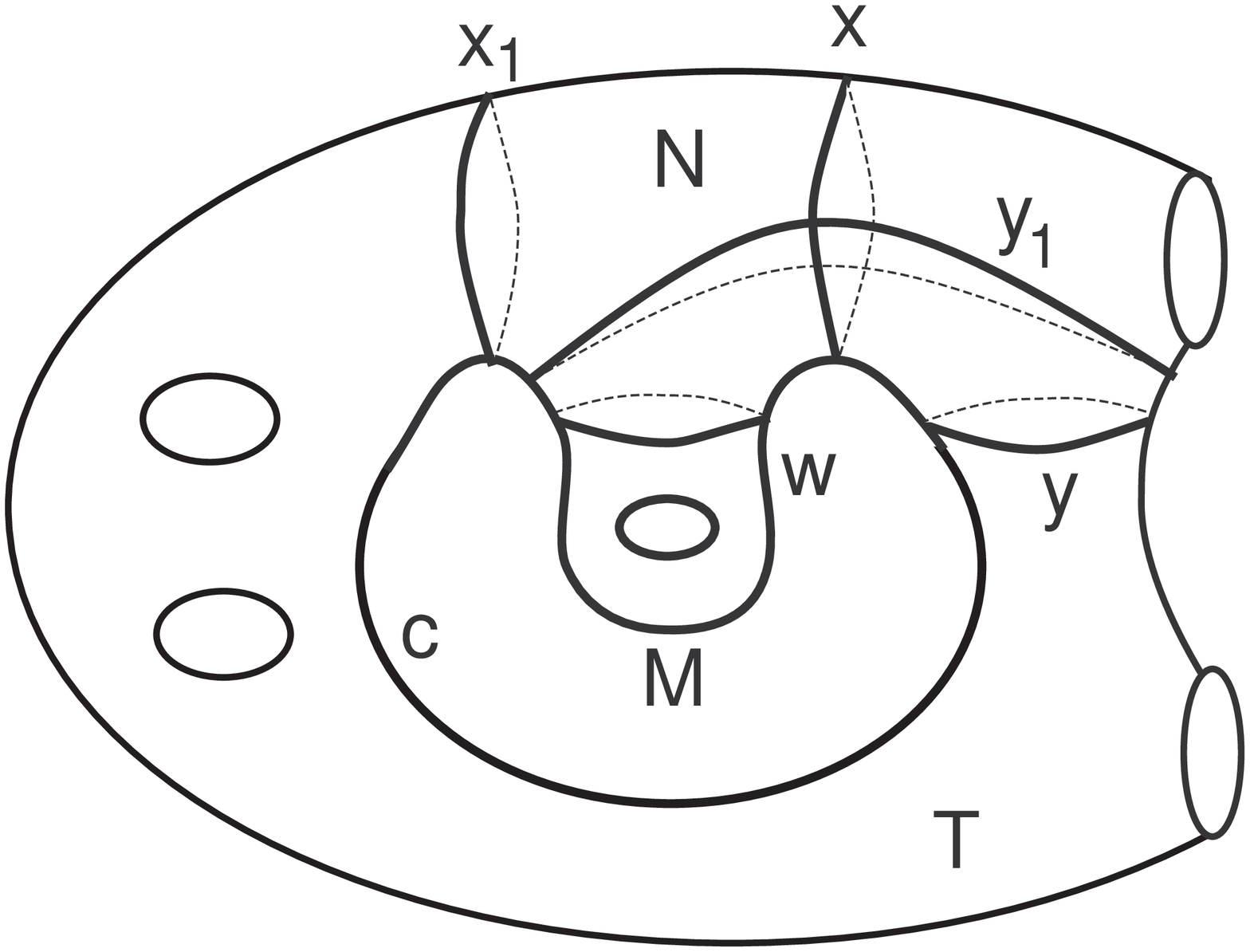} \caption{Curves
intersecting twice}
\end{center}
\end{figure}

In both cases we have seen that $(h_{x,y})_\# = (h_{x_1,y_1})_\#$.
By using our sequence, with an inductive argument we have $(h_{x,
y})_\# = (h_{z, t})_\#$, and $(h_{x, y})_\# = (h_{z, t})_\# =
\lambda$ on $\{[x]\} \cup \{[y]\} \cup \{[z]\} \cup \{[t]\} \cup
L_{x, y} \cup D_{x, y} \cup L_{z, t} \cup D_{z, t}$. In particular
we see that $(h_{x, y})_\#$ agrees with $\lambda$ on any
nonseparating curve $z$, hence $(h_{x, y})_\#$ agrees with
$\lambda$ on $\mathcal{N}(R)$. This proves the lemma.\end{proof}

\begin{theorem} \label{last} Suppose that the genus of $R$ is
at least two and $R$ has at most $g-1$ boundary components. A
simplicial map $\lambda : \mathcal{N}(R) \rightarrow
\mathcal{N}(R)$ is superinjective if and only if $\lambda$ is
induced by a homeomorphism of $R$.
\end{theorem}

\begin{proof} If $\lambda$ is induced by a homeomorphism of $R$,
then it preserves disjointness and nondisjointness and hence it is
superinjective. Assume that $\lambda$ is superinjective. In the
cases when $R$ is a closed surface or when $R$ has exactly one
boundary component, by Lemma \ref{extension} $\lambda$ extends to
a superinjective simplicial map $\lambda_*$ on $\mathcal{C}(R)$.
Then $\lambda_*$ is induced by a homeomorphism $h:R \rightarrow
R$, i.e. $\lambda_*(\alpha) = h_\#(\alpha)$ for each vertex
$\alpha$ in $\mathcal{C}(R)$, by the main results in \cite{Ir1}
and \cite{Ir2} and by Theorem \ref{theorem3}. Hence $\lambda$ is
induced by the
homeomorphism $h$.\\

Assume that $g \geq 3$ and $p=2$. Let $x$ and $y$ be disjoint
nonseparating circles such that $x, y$ and a boundary component,
$\partial_1$, of $R$ bound a pair of pants on $R$. Let $x', y'$ be
disjoint representatives of $\lambda([x]), \lambda([y])$
respectively. By using Lemma \ref{peripheral} and knowing that
$\lambda$ preserves disjointness and nondisjointness, it is easy
to see that $\lambda$ maps $\mathcal{N}(R_{x \cup y})$ to
$\mathcal{N}(R_{x' \cup y'})$. Since every essential separating
curve on $R_{x \cup y}$ separates $R$ into two subsurfaces each of
which has genus at least one, using chains on these subsurfaces we
could extend $\lambda$ to a superinjective simplicial map
$\lambda_{x, y} : \mathcal{C}(R_{x \cup y}) \rightarrow
\mathcal{C}(R_{x' \cup y'})$ as in case 1 in Lemma
\ref{extension}. Then by the main results in \cite{Ir2}, there
exists a homeomorphism $h : R_{x \cup y} \rightarrow R_{x' \cup
y'}$ such that $h(x) = x'$, $h(y) = y'$ and $\lambda_{x, y}$ is
induced by $h$. Then the proof of the theorem follows from Lemma
\ref{imp}.\\

\begin{figure}
\begin{center}
\epsfxsize=3in \epsfbox{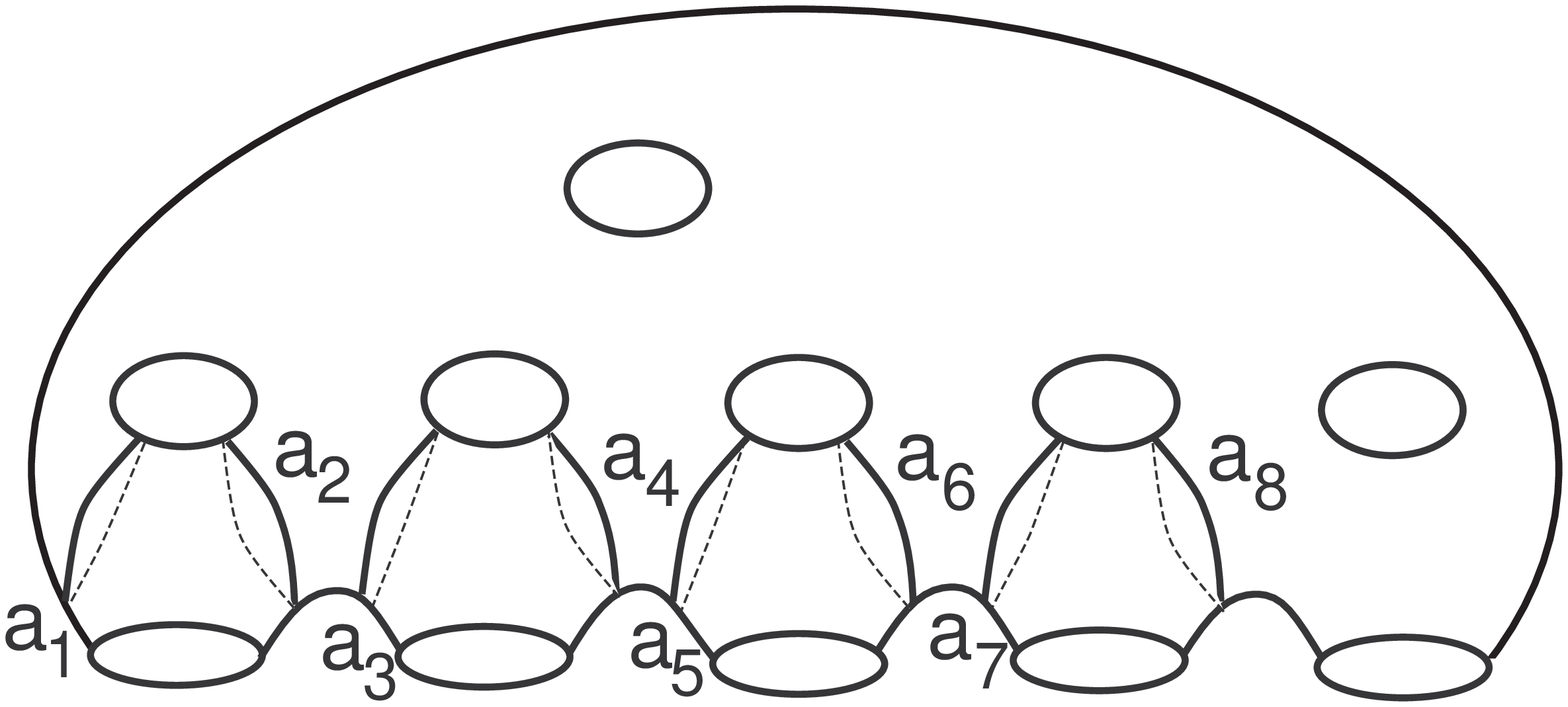} \caption{Cutting $R$ along
peripheral pairs}
\end{center}
\end{figure}

Now assume that $g \geq 4$ and $3 \leq p \leq g-1$. We will give
the proof when $p = g-1$. The proof of the remaining cases is
similar. Let $\{a_1, ..., a_{2(p-1)}\}$ be a set of pairwise
disjoint nonseparating circles such that $(a_{2i+1}, a_{2i+2})$ is
a peripheral pair as shown in Figure 13 for $i=0, ..., p-2$. Let
$R_{a_1 \cup a_2 ... \cup a_{2(p-1)}}$ be the genus two surface
with $2p-1$ boundary components which comes from the separation by
${a_1 \cup a_2 ... \cup a_{2(p-1)}}$. Exactly one of the boundary
components of $R_{a_1 \cup a_2 ... \cup a_{2(p-1)}}$ is a boundary
component of $R$. We identify $\mathcal{N}(R_{a_1 \cup a_2 ...
\cup a_{2(p-1)}})$ with a subcomplex $L_{a_1 \cup a_2 ... \cup
a_{2(p-1)}}$ of $\mathcal{N}(R)$. Let $\lambda_{a_1 \cup a_2 ...
\cup a_{2(p-1)}}$ denote the restriction of $\lambda$ on
$\mathcal{N}({R_{a_1 \cup a_2 ... \cup a_{2(p-1)}}})$. Let $a_1',
... , a_{2(p-1)}'$ be pairwise disjoint representatives of
$\lambda([a_1]),... ,$ $ \lambda([a_{2(p-1)}])$ respectively. By
using Lemma \ref{peripheral} and the properties of $\lambda$, it
is easy to see that $\lambda$ maps $\mathcal{N}({R_{a_1 \cup a_2
... \cup a_{2(p-1)}}})$ to $\mathcal{N}({R_{a_1' \cup a_2' ...
\cup a'_{2(p-1)}}})$. Since every essential separating curve on
$R_{a_1 \cup a_2 ... \cup a_{2(p-1)}}$ separates $R$ into two
subsurfaces each of which has genus at least one, using chains on
these subsurfaces we extend $\lambda$ to a superinjective
simplicial map $\lambda_{a_1 \cup a_2 ... \cup a_{2(p-1)}} :
\mathcal{C}(R_{a_1 \cup a_2 ... \cup a_{2(p-1)}}) \rightarrow
\mathcal{C}(R_{a_1' \cup a_2' ... \cup a'_{2(p-1)}})$ as in Lemma
\ref{extension}. Then by the main results in \cite{Ir2}, there
exists a homeomorphism $h : R_{a_1 \cup a_2 ... \cup a_{2(p-1)}}
\rightarrow R_{a_1' \cup a_2' ... \cup a'_{2(p-1)}}$ such that
$h(a_i) = a_i'$ and $\lambda_{a_1 \cup a_2 ... \cup a_{2(p-1)}}$
is induced by $h$. If we replace $(a_{2p-3}, a_{2p-2})$ with
another peripheral pair $(b_{2p-3}, b_{2p-2})$ where each of
$b_{2p-3}$ and $b_{2p-2}$ is disjoint from ${a_1 \cup a_2 ... \cup
a_{2p-4}}$, then by similar techniques we get a homeomorphism $t :
R_{a_1 \cup a_2 ... \cup a_{2p-4} \cup b_{2p-3} \cup b_{2p-2}}
\rightarrow R_{a_1' \cup a_2' ... \cup a'_{2p-4} \cup b'_{2p-3}
\cup b'_{2p-2}}$ such that $t(a_i) = a_i'$, $t(b_j) = b_j'$, where
$b'_{2p-3}, b'_{2p-2}$ are pairwise disjoint representatives of
$\lambda([b_{2p-3}]), \lambda([b_{2p-2}])$ respectively, and the
map $\lambda_{a_1 \cup a_2 ... \cup a_{2p-4} \cup b_{2p-3} \cup
b_{2p-2}}$ is induced by $t$. Then by following the techniques in
the proof of Lemma \ref{imp}, we see that there exists a
homeomorphism $r : R_{a_1 \cup a_2 ... \cup a_{2p-4}} \rightarrow
R_{a_1' \cup a_2' ... \cup a'_{2p-4}}$ such that $r(a_i) = a_i'$
and $\lambda_{a_1 \cup a_2 ... \cup a_{2p-4}}$ is induced by $r$.
Then by an inductive argument, there exists a homeomorphism $q : R
\rightarrow R$ such that $\lambda$ is induced by $q$.\end{proof}\\

Now we consider the graph $\mathcal{G}(R)$ defined by Schaller.
The vertex set of $\mathcal{G}(R)$ is the set of isotopy classes
of nontrivial nonseparating simple closed curves on $R$. Two
vertices are connected by an edge if and only if their geometric
intersection number is one. In the proof of the following theorem,
we use some of Schaller's results given in \cite{Sc} that if $g
\geq 2$ and $R$ is not a closed surface of genus two, then
$Aut(\mathcal{G}(R))= Mod_R^*$.

\begin{theorem} \label{lastlast} Suppose that $R$ has genus at least two.
If $R$ is not a closed surface of genus two, then
$Aut(\mathcal{N}(R))= Mod_R^*$. If $R$ is a closed surface of
genus two, then $Aut(\mathcal{N}(R))= Mod_R ^* /\mathcal{C}(Mod_R
^*)$.\end{theorem}

\begin{proof} Assume that $R$ is not a closed surface of genus
two. If $[f] \in Mod_R^*$, then $[f]$ induces an automorphism of
$\mathcal{N}(R)$. If $[f]$ fixes the isotopy class of every
nontrivial nonseparating simple closed curve on $R$, then $f$ is
orientation preserving and it can be shown that it is isotopic to
$id_R$. An automorphism $\lambda$ of $\mathcal{N}(R)$ is a
superinjective simplicial map of $\mathcal{N}(R)$. By Lemma
\ref{intone}, it preserves geometric intersection one property,
and hence induces an automorphism of $\mathcal{G}(R)$. An
automorphism of $\mathcal{G}(R)$ is induced by a homeomorphism of
$R$, \cite{Sc} (Note that if $R$ has at most $g-1$ boundary
components, then Theorem \ref{last} also implies that $\lambda$ is
induced by a homeomorphism of $R$). Then it is easy to see that we
get $Aut(\mathcal{N}(R))= Mod_R^*$. Assume that $R$ is a closed
surface of genus two. If $[f] \in Mod_R^*$, then $[f]$ induces an
automorphism of $\mathcal{N}(R)$. If $[f]$ fixes the isotopy class
of every nonseparating simple closed curve on $R$ then $f$ is
either isotopic to $id_R$ or the hyperelliptic involution. An
automorphism of $\mathcal{N}(R)$ is a superinjective simplicial
map of $\mathcal{N}(R)$ and by Theorem \ref{last} it is induced by
a homeomorphism of $R$. Then we have $Aut(\mathcal{N}(R))= Mod_R^*
/\mathcal{C}(Mod_R ^*)$.\end{proof}\\

Note that by the main results of this paper and the results in
\cite{Ir1}, \cite{Ir2}, \cite{Iv1}, \cite{Sc}, we have that if $R$
has genus at least two, then $Aut(\mathcal{N}(R))=
Aut(\mathcal{C}(R)) =
Aut(\mathcal{G}(R))$.\\

{\bf Acknowledgments}\\

We thank John D. McCarthy for his suggestions and many discussions
about this work. We also thank Peter Scott for his comments.

\vspace{0.2cm}

\noindent University of Michigan, Department of Mathematics, Ann
Arbor, MI 48109, USA;

\noindent eirmak@umich.edu\\

\end{document}